\documentclass[12pt,reqno]{amsart}

\usepackage{amsfonts,color,amsthm,amsmath,amssymb}

\usepackage{color}
\def\Box{\vcenter{\vbox{\hrule\hbox{\vrule
     \vbox to 8.8pt{\hbox to 10pt{}\vfill}\vrule}\hrule}}}

\newcommand{\tr}{\text{Tr}}
\def\qed{{\hfill$\square$}}
\def\proof{{\vspace{-0.0cm}\bf Proof: \,}}

\def\Z{{\mathbb Z}}

\def\C{{\mathbb C}}
\def\F{{\mathbb F}}

\def\mod{{\mathrm{mod\,\,}}}

\def\Tr{{\mathrm{Tr}}}
\def\Norm{{\mathrm{Norm}}}
\def\Cay{{\mathrm{Cay}}}

\def\PG{{\mathrm{PG}}}

\def\AP{{\mathrm{AP}}}

\newcommand\cQ{{\mathcal Q}}

\newtheorem{theorem}{Theorem}[section]

\newtheorem{lemma}[theorem]{Lemma}
\newtheorem{remark}[theorem]{Remark}

\newtheorem{corollary}[theorem]{Corollary}

\newtheorem{proposition}[theorem]{Proposition}

\newtheorem{conj}[theorem]{Conjecture}
\numberwithin{equation}{section}

\begin{document}

\title[Strongly regular Cayley graphs]{Strongly regular Cayley graphs from partitions of subdifference sets of the Singer difference sets}

\author[Momihara and Xiang]{Koji Momihara$^*$, Qing Xiang$^{\dagger}$}

\address{Koji Momihara,  Faculty of Education, Kumamoto University, 2-40-1 Kurokami, Kumamoto 860-8555, Japan} \email{momihara@educ.kumamoto-u.ac.jp}

\address{Qing Xiang, Department of Mathematical Sciences, University of Delaware, Newark, DE 19716, USA} \email{qxiang@udel.edu}


\thanks{$^\ast$Research supported by JSPS under Grant-in-Aid for Young Scientists (B) 17K14236 and Scientific Research (B) 15H03636.} 

\thanks{$^{\dagger}$Research partially supported by an NSF grant DMS-1600850}


\begin{abstract}
In this paper, we give a new lifting construction of ``hyperbolic" type of strongly regular Cayley graphs. Also we give new constructions of strongly regular Cayley graphs over the additive groups of finite fields based on partitions of subdifference sets of the Singer difference sets. Our results unify some recent constructions of strongly regular Cayley graphs related to $m$-ovoids and $i$-tight sets in finite geometry. Furthermore,  some of the strongly regular Cayley graphs obtained in this paper are new or nonisomorphic to known strongly regular graphs with the same parameters.  
\end{abstract}

\maketitle

\section{Introduction}\label{sec:1}
We assume that the reader is familiar with the basic theory of strongly regular graphs and difference sets. For strongly regular graphs (srgs), our main references are \cite{BH} and \cite{cg}. For difference sets, we refer the reader to \cite{Lander} and Chapter 6 of \cite{bjl}. Strongly regular graphs are closely related to many other combinatorial/geometric objects, such as two-weight codes, two-intersection sets, $m$-ovoids, $i$-tight sets, and partial difference sets. For these connections, we refer the reader to \cite[p.~132]{BH}, \cite{CK86, Ma}, and some more recent papers \cite{FMX, DDMR, BLMX} on Cameron-Liebler line classes and hemisystems.

Let $\Gamma$ be a (simple, undirected) graph. The adjacency matrix of $\Gamma$ is the $(0,1)$-matrix $A$ with both rows and columns indexed by the vertex set of $\Gamma$, where $A_{xy} = 1$ when there is an edge between $x$ and $y$ in $\Gamma$ and $A_{xy} = 0$ otherwise. A useful way to check whether a graph is strongly regular is by using the eigenvalues of its adjacency matrix. For convenience we call an eigenvalue {\it restricted} if it has an eigenvector perpendicular to the all-ones vector ${\bf 1}$. (For a $k$-regular connected graph, the restricted eigenvalues are the eigenvalues different from $k$.)

\begin{theorem}\label{char}
For a simple graph $\Gamma$ of order $v$, not complete or edgeless, with adjacency matrix $A$, the following are equivalent:
\begin{enumerate}
\item $\Gamma$ is strongly regular with parameters $(v, k, \lambda, \mu)$ for certain integers $k, \lambda, \mu$,
\item $A^2 =(\lambda-\mu)A+(k-\mu) I+\mu J$ for certain real numbers $k,\lambda, \mu$, where $I, J$ are the identity matrix and the all-ones matrix, respectively, 
\item $A$ has precisely two distinct restricted eigenvalues.
\end{enumerate}
\end{theorem}

For a proof of Theorem~\ref{char}, we refer the reader to \cite{BH}. An effective method to construct strongly regular graphs is by using Cayley graphs. Let $G$ be an additively written group of order $v$, and let $D$ be a subset of $G$ such that $0\not\in D$ and $-D=D$, where $-D=\{-d\mid d\in D\}$. The {\it Cayley graph over $G$ with connection set $D$}, denoted ${\rm Cay}(G,D)$, is the graph with the elements of $G$ as vertices; two vertices are adjacent if and only if their difference belongs to $D$. In the case when $\Cay(G,D)$ is a strongly regular graph, the connection set $D$ is called a (regular) {\it partial difference set}. Examples of strongly regular Cayley graphs are the Paley graphs ${\rm P}(q)$, where $q$ is a prime power congruent to 1 modulo 4, the Clebsch graph, and the affine orthogonal graphs (\cite{BH}). For $\Gamma={\rm Cay}(G,D)$ with $G$ abelian, the eigenvalues of $\Gamma$ are exactly $\chi(D):=\sum_{d\in D}\chi(d)$, where $\chi$ runs through the character group of  $G$. This fact reduces the problem of computing eigenvalues of abelian Cayley graphs to that of computing some character sums, and is the underlying reason why the Cayley graph construction has been very effective for the purpose of constructing srgs. The survey of Ma~\cite{Ma} contains much of what is known about partial difference sets and about connections with strongly regular graphs. 

A $(v,k,\lambda,\mu)$ srg is said to be of
{\em{Latin square type}} (respectively, {\it negative Latin square
type}) if $(v,k,\lambda,\mu) = (n^2, r(n-\epsilon), \epsilon n+r^2-3
\epsilon r, r^2 - \epsilon r)$ and $\epsilon = 1$ (respectively, $\epsilon=-1$). When $v$ (the number of vertices) is a prime power, many constructions of srgs with Latin square or negative Latin square type parameters are known. For example, the srgs arising from partial spreads of $\PG(2m -1,q)$ have Latin square parameters, and the affine orthogonal graphs, ${\rm VO}^{-}(2m,q)$, have negative Latin square type parameters. Still the range of $r$ in the parameters $(n^2, r(n-\epsilon), \epsilon n+r^2-3
\epsilon r, r^2 - \epsilon r)$ of the known srgs of Latin square or negative Latin square type can sometimes be limited; moreover Latin square and negative Latin square type strongly regular Cayley graphs with certain extra properties\footnote{For example, the elements of the connection set must all lie on a quadratic surface.} have found many connections with finite geometric objects such as $m$-ovoids and $i$-tight sets (cf. \cite{FMX, DDMR, BLMX}). Therefore it is of interest to construct more strongly regular Cayley graphs of Latin square or negative Latin square type. The purpose of the current paper is two fold. First, we give new constructions of strongly regular Cayley graphs, and obtain some new srgs. Secondly, we unify and give simpler proofs for some recent constructions of strongly regular Cayley graphs. 

The paper is organized as follows. In Section 2, we review some basic properties of Gauss sums which will be used in later sections. In Section 3, we give two constructions of strongly regular Cayley graphs on the additive group of $\F_{q^2}$ by lifting a cyclotomic strongly regular graph on $\F_q$. The first lifting construction (Proposition~\ref{theorem:main1}) is of ``elliptic" type, and it was already given in \cite{MX}. The second lifting construction (Proposition~\ref{theorem:main2}) is of ``hyperbolic" type, and this constrution is new. In Section 4, we 
generalize and unify the constructions of strongly regular Cayley graphs corresponding to $m$-ovoids and $i$-tight sets in \cite{FMX, BLMX}. We give a general construction of strongly regular Cayley graphs by using a certain partition of a subdifference set (and its complement) of the Singer difference set. When the subdifference sets arise from subfields, we recover the results in \cite{FMX, BLMX}. In Sections 5, 6, and 7, we apply the general construction in Section 4 to the three known cases of subdifference sets of the Singer difference sets, namely, the semiprimitive case, the sporadic case, and the subfield case. We either recover strongly regular Cayley graphs constructed in some of our recent papers \cite{Mo, FMX, BLMX}, or we produce new strongly regular Cayley graphs. In particular, Corollaries 7.7 and 7.9 give strongly regular Cayley graphs with the same parameters as the affine polar graphs. By using a computer, it is shown that the newly constructed graphs in Corollaries 7.7 and 7.9 are not isomorphic to the affine polar graphs when the parameters are small.

\section{Preliminaries}
We will use Gauss sums and Gauss periods to compute character values of certain subsets of  $\F_{q}$, the finite field of order $q$. So it is helpful to introduce characters of both kinds of finite fields, and review basic properties of Gauss sums. Let $p$ be a prime, $f$ a positive integer, and $q=p^f$. The canonical additive character $\psi_{\F_q}$ of $\F_q$ is defined by 
$$\psi_{\F_q}\colon\F_q\to \C^{\ast},\qquad\psi_{\F_q}(x)=\zeta_p^{\tr _{q/p}(x)},$$
where $\zeta_p={\rm exp}(\frac {2\pi i}{p})$ is a complex primitive $p$-th root of unity and $\Tr _{q/p}$ is the trace function from $\F_q$ to $\F_p$
defined by 
$\tr_{q/p}(x) = x + x^p + x^{p^2} +\cdots + x^{p^{f-1}}$. 
All the additive characters of $\F_q$ can be obtained from the canonical one. For $a\in \F_q$, define 
\begin{equation}\label{additive}
\psi_a(x)=\psi_{\F_q}(ax), \;\forall x\in \F_q.
\end{equation}
Then $\{\psi_a\mid a\in \F_q\}$ is the group of additive characters of $\F_q$. For a multiplicative character
$\chi$ and the canonical
additive character $\psi_{\F_q}$ of $\F_q$, define the {\it Gauss sum} by
\[
G_q(\chi)=\sum_{x\in \F_q^\ast}\chi(x)\psi_{\F_q}(x).
\]

Some basic properties of Gauss sums are listed below:
\begin{proposition}{\em (\cite[Theorem 5.2]{LN97})}
Let $\chi$ be a multiplicative character of $\F_{q}$. Then, the following hold:  
\begin{enumerate}
\item[(i)] $G_q(\chi)\overline{G_q(\chi)}=q$ if $\chi$ is nontrivial;
\item[(ii)] $G_q(\chi^p)=G_q(\chi)$, where $p$ is the characteristic of $\F_q$;
\item[(iii)] $G_q(\chi^{-1})=\chi(-1)\overline{G_q(\chi)}$;
\item[(iv)] $G_q(\chi)=-1$ if $\chi$ is trivial.
\end{enumerate}
\end{proposition}
Let $\omega$ be a fixed primitive element of $\F_q$ and $N$ a positive integer dividing $q-1$. For $0\le i\le N-1$ we set $C_i^{(N,q)}=\omega^i C_0$, where $C_0$ is the subgroup of index $N$ of $\F_q^\ast$. The {\it Gauss periods} associated with these cosets are defined by $\psi_{\F_q}(C_i^{(N,q)}):=\sum_{x\in C_i^{(N,q)}}\psi_{\F_q}(x)$, $0\le i\le N-1$, where $\psi_{\F_q}$ is the canonical additive character of $\F_q$. By orthogonality of characters, the Gauss periods can be expressed as a linear combination of Gauss sums:
\begin{equation}
\psi_{\F_q}(C_i^{(N,q)})=\frac{1}{N}\sum_{j=0}^{N-1}G_q(\chi^{j})\chi^{-j}(\omega^i), \; 0\le i\le N-1,
\end{equation}
where $\chi$ is any fixed multiplicative character of order $N$ of $\F_q$.  
For example, if $N=2$, 
we have
\begin{equation}\label{eq:Gaussquad}
\psi_{\F_q}(C_i^{(2,q)})=\frac{-1+(-1)^iG_q(\eta)}{2},\; 0\le i\le 1,
\end{equation}
where $\eta$ is the quadratic character of $\F_q$. 

The quadratic Gauss sum, $G_q(\eta)$, can be evaluated explicitly. 
\begin{theorem}\cite[Theorem~5.15]{LN97} \label{thm:Gauss}
Let $q=p^s$ be a prime power with $p$ a prime and $\eta$ be the quadratic character of $\F_q$. 
Then, 
\begin{equation}\label{eq:Gaussquad1}
G_q(\eta)=\begin{cases}
(-1)^{s-1}q^{1/2}& \text{ if }  p\equiv 1\,(\mod{4}), \\
(-1)^{s-1}i^s q^{1/2} & \text{ if } p\equiv 3\,(\mod{4}). 
\end{cases}
\end{equation}
\end{theorem}
Also, in the semi-primitive case, the 
Gauss sum can be computed.
\begin{theorem}\label{thm:semiprim}{\em (\cite[Theorem~11.6.3]{BEW97})}
Let $p$ be a prime. 
Suppose that $N>2$ and $p$ is semi-primitive modulo $N$, 
i.e., there exists a positive integer $j$ such that  $p^j\equiv -1\pmod{N}$. Choose 
$j$ minimal and write 
$f=2js$ for any positive integer $s$. Let $\chi$ be a multiplicative character of order $N$ of $\F_{p^f}$. 
Then, 
\[
p^{-f/2}G_{p^f}(\chi)=
\left\{
\begin{array}{ll}
(-1)^{s-1},&  \mbox{if $p=2$,}\\
(-1)^{s-1+(p^j+1)s/N},&  \mbox{if $p>2$. }
 \end{array}
\right.
\]
\end{theorem}

The following theorems are referred to as the  {\it Davenport-Hasse lifting formula} and the  {\it Davenport-Hasse product formula}, respectively. 
\begin{theorem}\label{thm:lift}{\em (\cite[Theorem~5.14]{LN97})}
Let $m$ be a positive integer. Let $\chi$ be a nontrivial multiplicative character of $\F_{q}$ and $\chi'$ 
be the lift of $\chi$ to $\F_{q^m}$, i.e., $\chi'(x)=\chi(\Norm_{\F_{q^m}/\F_q}(x))$
for $x\in \F_{q^m}$, where $\Norm_{\F_{q^m}/\F_q}$ is the norm from $\F_{q^m}$ to $\F_{q}$. Then, 
\[
G_{q^m}(\chi')=(-1)^{m-1} G_{q}(\chi)^m. 
\]
\end{theorem} 
\begin{theorem}
\label{thm:Stickel2}{\em (\cite[Theorem~11.3.5]{BEW97})}
Let $\eta$ be a multiplicative character of order $\ell>1$ of  $\F_{q}$. For  every nontrivial multiplicative character $\chi$ of $\F_{q}$, 
\[
G_{q}(\chi)=\frac{G_{q}(\chi^\ell)}{\chi^\ell(\ell)}
\prod_{i=1}^{\ell-1}
\frac{G_{q}(\eta^i)}{G_{q}(\chi\eta^i)}. 
\]
\end{theorem}

We will use the following formula later. 
\begin{theorem}\label{prop:charaadd}{\em (\cite[Theorem~5.33]{LN97})}
Let $\psi_{\F_q}$ be the canonical additive character of $\F_q$ with $q$ odd, and 
let $f(x)=a_2x^2+a_1x+a_0\in \F_q[x]$ with $a_2\not=0$. Then
\[
\sum_{x\in \F_q}\psi_{\F_q}(f(x))=\psi_{\F_q}(a_0-a_1^2(4a_2)^{-1})\eta(a_2)G_q(\eta), 
\]
where $\eta$ is the quadratic character of $\F_q$. 
\end{theorem}
\section{Basic lifting constructions}
\subsection{Subdifference sets of the Singer difference sets}\label{sec:Quo}
Let $p$ be a prime, $f\geq 1$, $m\geq 2$ be integers and $q=p^f$. Let $L$ be a complete system of coset representatives of $\F_q^\ast$ in $\F_{q^m}^\ast$. We can, and do, choose $L$ in such a way that $\Tr_{q^m/q}(x)=0$ or $1$ for any $x\in L$. 
Let 
\[
L_0=\{x\in L\,|\,\Tr_{q^m/q}(x)=0\}\mbox{ and } L_1=\{x\in L\,|\,\Tr_{q^m/q}(x)=1\}. 
\]
Then, 
\begin {equation*}\label{Singer}
H_0=\{\overline{x}\in \F_{q^m}^\ast/\F_q^\ast\,|\,x\in L_0\}
\end{equation*}
represents a hyperplane of the projective space $\PG(m-1,q)$, and it is a so-called {\it Singer difference set} in the cyclic group $\F_{q^m}^\ast/\F_q^\ast$. (Here $\overline{x}=x\F_q^\ast$ represents the 
projective point corresponding to the one-dimensional subspace over $\F_q$ spanned by $x$.)  

Any nontrivial multiplicative character $\chi$ of exponent $(q^m-1)/(q-1)$ of $\F_{q^m}^\ast$ induces a character of the quotient group $\F_{q^m}^\ast/\F_q^\ast$, which will also denoted by $\chi$. Moreover, every character of $\F_{q^m}^\ast/\F_q^\ast$ arises in this way. By a result given in \cite{Y}, for any nontrivial multiplicative character $\chi$ of exponent $(q^m-1)/(q-1)$ of $\F_{q^m}^\ast$, we have 
$$\chi(H_0)=G_{q^m}(\chi)/q.$$

Assume that $N\,|\,(q^m-1)/(q-1)$. Then $$\overline{C_0}:=C_0^{(N,q^m)}/\F_q^\ast\le \F_{q^m}^\ast/\F_q^\ast.$$ 
Let $S$ be a complete system of coset representatives of $\overline{C_0}$ in $\F_{q^m}^\ast/\F_q^\ast$,  and set $G=\{s\overline{C_0}\,|\,s\in S\}\simeq \F_{q^m}^\ast/C_0^{(N,q^m)}$. For convenience, we will use $\tilde{s}$ to denote  
$s\overline{C_0}$. 

In the rest of this section, we will assume that $\Cay(\F_{q^m},C_0^{(N,q^m)}) $ is strongly regular, where $N\,|\,(q^m-1)/(q-1)$. Such a strongly regular graph is called {\it cyclotomic}. The following three series of cyclotomic strongly regular graphs were known~\cite{SW}: 
\begin{itemize}
\item[(1)] (subfield case) $C_0^{(N,q^m)}=\F_{q^d}^\ast$ where $d\,|\,m$,  
\item[(2)] (semi-primitive case) $-1\in \langle p\rangle\le (\Z/N\Z)^\ast$,
\item[(3)] (sporadic case) $\Cay(\F_{q^m},C_0^{(N,q^m)})$ has one of the eleven sets of parameters 
given in Table~\ref{Tab1}. 
\begin{table}[h]
\caption{Eleven sporadic examples}
\vspace{-0.5cm}
\label{Tab1}
$$
\begin{array}{|c||c|c|c|c|c|c|c|c|c|c|c|}
\hline
N&     11&19 &35      &37 &43   &67       &107&133&163&323&499\\
\hline
q^m&3^5&5^9&3^{12}&7^9&11^7&17^{33}&3^{53}&5^{18}&41^{81}&3^{144}
&5^{249}\\
\hline
\end{array}
$$
\end{table}
\end{itemize}
We mention in passing that Schmidt and White~\cite{SW} conjectured that besides the above three cases, there are no more cyclotomic strongly regular graphs. 
\begin{conj}
Let $N\,|\,\frac{q^m-1}{q-1}$ with $N>1$. If $\Cay(\F_{q^m},C_0^{(N,q^m)})$ is 
strongly regular, then one of (1), (2), or (3) above holds. 
\end{conj}

This conjecture remains open. Some partial results were obtained in \cite{SW}.

Assume that $\Cay(\F_{q^m},C_0^{(N,q^m)}) $ is strongly regular, where $N\,|\,\frac{q^m-1}{q-1}$. Then $|H_0\cap s\overline{C_0}|$, $s\in S$, take exactly two values. It follows that $|H_0\cap s\overline{C_0}|-|H_0\cap \overline{C_0}|=0$ or $\delta$, where $\delta$ is a nonzero integer. For any nontrivial multiplicative character $\chi$ of $\F_{q^m}$ of exponent $N$, 
\begin{align*}
\chi(H_0)=\,&\sum_{s\in S}|H_0\cap s\overline{C_0}|\chi(\tilde{s})\\
=\,&\sum_{s\in S}(|H_0\cap s\overline{C_0}|-|H_0\cap \overline{C_0}|)\chi(\tilde{s})\\
=\,&\delta\sum_{s\in S'}\chi(\tilde{s}), 
\end{align*}
where 
\begin{equation}\label{def:S'}
S'=\{s\in S : |H_0\cap s\overline{C_0}|-|H_0\cap \overline{C_0}|=\delta\}.
\end{equation}
 Thus
\begin{equation}\label{eq:sum}
\sum_{s\in S'}\chi(\tilde{s})=\frac{\chi(H_0)}{\delta}=\frac{G_{q^m}(\chi)}{\delta q}. 
\end{equation}
It follows that $\delta$ is a power of $p$. Furthermore, noting that $G_{q^m}(\chi)\overline{G_{q^m}(\chi)}=q^m$, we see that the set $\{\tilde{s}\mid s\in S'\}\subset G$ forms a difference set, 
which is called a {\it subdifference set} of $H_0$~\cite{Mc}. Let $\omega$ be a primitive element of $\F_{q^m}$. Then we could choose $S=\{\overline{1}, \overline{\omega}, \ldots , \overline{\omega}^{N-1}\}$, where $\overline{\omega}=\omega\F_q^\ast$. In this way, since $S'$ is a subset of $S$, we define  
\begin{equation}\label{eq:subdi}
I=\{0\le i\le N-1\mid \overline{\omega}^i\in S'\}.
\end{equation}
In the rest of the paper, we will also call $I$ a subdifferecne set in $\Z_N$ of the Singer difference set. 
\subsection{Two lifting constructions}\label{sec:two}
Let $\gamma$ be a primitive element of $\F_{q^{2m}}$ and set $\omega=\gamma^{q^m+1}$. Then, $\omega$ is a primitive element in $\F_{q^m}$. Let 
$C_j^{(N,q^{2m})}=\gamma^j\langle \gamma^N\rangle$, $0\le j\le N-1$. 
The following lifting construction was already given in~\cite{MX}. For completeness, we repeat the construction here.
\begin{proposition} \label{theorem:main1}
Assume that $\F_q^\ast\le C_0^{(N,q^m)}\le \F_{q^m}^\ast$ and $\Cay(\F_{q^m},C_0^{(N,q^m)})$ is strongly regular. Let $I$ be the corresponding subdifferecne set defined in \eqref{eq:subdi}. 
Let 
\begin{equation}\label{ellip}
E=\bigcup_{i\in I}C_i^{(N,q^{2m})}. 
\end{equation}
Then $\Cay(\F_{q^{2m}}, E)$ is a strongly regular graph with negative Latin square type parameters $(n^2,r(n+1),-n+r^2+3r,r^2+r)$, where $n=q^m$ and $r=(q^m-1)|I|/N$. 
\end{proposition}
\proof
Let $\psi_{\F_{q^{2m}}}$ be the canonical additive character of $\F_{q^{2m}}$ and let  $\chi_N'$ be a multiplicative character of order $N$ of $\F_{q^{2m}}$. We will show that $\psi_{\F_{q^{2m}}}(\gamma^a E)$, $0\leq a\leq N-1$, take exactly two distinct values. By the orthogonality of characters, we compute 
\begin{eqnarray*}
S_a=N\cdot \psi_{\F_{q^{2m}}}(\gamma^a E)+|I|=\sum_{j=1}^{N-1}G_{q^{2m}}(\chi_N'^{-j})\sum_{i\in I}\chi_N'^{j}(\gamma^{a+i}). 
\end{eqnarray*}
Since $N\,|\,\frac{q^m-1}{q-1}$, there is a multiplicative character $\chi_N$ of $\F_{q^m}$ of order $N$ such that $\chi_N'(\gamma)=\chi_N(\omega)$, i.e., $\chi_N'$ is the lift of $\chi_N$. 
By the Davenport-Hasse lifting formula, we have 
\[
S_a
=-\sum_{j=1}^{N-1}{\chi_{N}}^{j}(\omega^{a})G_{q^m}(\chi_{N}^{-j})G_{q^m}(\chi_{N}^{-j})
\sum_{i\in I}\chi_{N}^{j}(\omega^{i}). 
\]
On the other hand, from the definition of $I$, we have 
\begin{equation}\label{eq:defIuse}
\sum_{i\in I}\chi_N^j(\omega^i)=\sum_{s\in S'}\chi_N^j(\tilde{s})=\frac{G_{q^m}(\chi_N^j)}{\delta q}. 
\end{equation}
Hence, 
\begin{align*}
S_a
=\,&-\frac{1}{\delta q} \sum_{j=1}^{N-1}\chi_{N}^{j}(\omega^{a})G_{q^m}(\chi_{N}^{-j})G_{q^m}(\chi_{N}^{-j})G_{q^m}(\chi_{N}^{j})\nonumber\\
=\,&- \frac{q^{m-1}}{\delta} \sum_{j=1}^{N-1}\chi_{N}^{j}(\omega^{a})G_{q^m}(\chi_{N}^{-j})\\
=\,&-q^m\sum_{j=1}^{N-1}\sum_{i \in I}\chi_N^{-j}(\omega^i)\chi_N^{j}(\omega^a)
=q^m|I|-\begin{cases}
q^m N,& \text{ if }  a\in I, \\
0,& \text{ if } a\not\in I. 
\end{cases}  
\end{align*}
Thus,  $\psi_{\F_{q^{2m}}}(\gamma^a E)=\frac{S_{a}-|I|}{N}$, $0\leq a\leq N-1$, take exactly two distinct values 
$\frac{(q^m-1)|I|}{N}$ and $\frac{(q^m-1)|I|}{N}-q^m$. 
Therefore $\Cay(\F_{q^{2m}}, E)$ is strongly regular. The parameters of $\Cay(\F_{q^{2m}}, E)$ can be computed in a straightforward way. We omit the details. \qed
\begin{remark}{\em 
If $\Cay(\F_{q^m}, C_0^{(N,q^m)})$ is a cyclotomic strongly regular graph in the subfield case with $N=\frac{q^m-1}{q-1}$, then $C_0^{(N,q^m)}=\F_q^\ast$, $S=\F_{q^m}^\ast/\F_q^\ast$, and $S'=H_0$. In this case, we find that
\[
E=\{x\in \F_{q^{2m}}^\ast\,|\,\Tr_{q^m/q}(x^{q^m+1})=0\},
\] where $\Tr_{q^m/q}(x^{q^m+1})$ is a nondegenerate $\F_q$-valued elliptic quadratic form on $\F_{q^{2m}}$. Therefore it is appropriate to call the lifting construction given in Proposition~\ref{theorem:main1} an {\it elliptic} type lifting construction.}
\end{remark}
We give a new lifting construction, which is  of ``hyperbolic" type. 
\begin{proposition} \label{theorem:main2}
Let $\omega$ be a primitive element of $\F_{q^m}$. Assume that $\F_q^\ast\le C_0^{(N,q^m)}\le \F_{q^m}^\ast$ and $\Cay(\F_{q^m},C_0^{(N,q^m)})$ is strongly regular. Let $I$ be the corresponding subdifferecne set defined in \eqref{eq:subdi}. Let
\begin{equation}\label{hyperb}
H=\{(y,y^{-1}x\omega^\ell)\,|\,x \in C_0^{(N,q^m)},y \in \F_{q^m}^\ast,\ell \in I\}\subseteq \F_{q^m}\times \F_{q^m}. 
\end{equation}
Then $\Cay(\F_{q^m} \times \F_{q^m}, H)$ is a strongly regular graph with Latin square type parameters  $(n^2,r(n-1),n+r^2-3r,r^2-r)$, where $n=q^m$ and $r=(q^m-1)|I|/N$. 
\end{proposition}
\proof
Let $\psi_{\F_{q^m}}$ be the canonical additive character of $\F_{q^{m}}$ and let  $\chi_{q^m-1}$ be a multiplicative character of order $q^m-1$ of $\F_{q^{m}}$. 
Each additive character of $\F_{q^m}\times \F_{q^m}$ has the form
\begin{equation}\label{fieldchara22}
\psi_{a,b}((x,y))=\psi_{\F_{q^m}}(ax+by),\; \quad (x,y)\in \F_{q^m}\times \F_{q^m},
\end{equation}
where $(a,b)\in \F_{q^m}\times \F_{q^m}$. 
Then, by the definition of $H$, we need to compute the character values: 
\[
S_{a,b}:=\sum_{y\in \F_{q^m}^\ast}\sum_{x\in C_0^{(N,q^m)}}\sum_{\ell\in I}\psi_{\F_{q^m}}(ay+bxy^{-1}\omega^\ell), \, \, \quad (0,0)\neq (a,b)\in \F_{q^m}\times \F_{q^m}.
\]
In the case where either one of $a$ or $b$ is zero, it is clear that 
$S_{a,b}=-(q^m-1)|I|/N$. 

Now, we assume that $a\not=0$ and $b\not=0$. 
By the orthogonality of characters, we have 
\begin{equation}\label{eq:hypfir1}
S_{a,b}=\frac{1}{(q^m-1)^2}\sum_{j,k=0}^{q^m-2}
\sum_{y\in \F_{q^m}^\ast}\sum_{x\in C_0^{(N,q^m)}}\sum_{\ell\in I}G_{q^m}(\chi_{q^m-1}^{-j})G_{q^m}(\chi_{q^m-1}^{-k})\chi_{q^m-1}^{j}(a)
\chi_{q^m-1}^{k}(bx\omega^\ell)\chi_{q^m-1}^{j-k}(y).
\end{equation}
Since $\sum_{y\in \F_{q^m}^\ast}\chi_{q^m-1}^{j-k}(y)=q^m-1$ or $0$ according as $j\equiv k\,(\mod{q^m-1})$ or $j\not\equiv k\,(\mod{q^m-1})$, continuing from \eqref{eq:hypfir1}, 
we have 
\begin{equation}\label{eq:hypfir2}
S_{a,b}=\frac{1}{q^m-1}\sum_{j=0}^{q^m-2}
\sum_{x\in C_0^{(N,q^m)}}\sum_{\ell\in I}G_{q^m}(\chi_{q^m-1}^{-j})^2\chi_{q^m-1}^{j}(a)
\chi_{q^m-1}^{j}(bx\omega^\ell). 
\end{equation}
Let $\chi_{N}:=\chi_{q^m-1}^{\frac{q^m-1}{N}}$. 
Since $\sum_{x\in C_0^{(N,q^m)}}\chi_{q^m-1}^{j}(x)=(q^m-1)/N$ or $0$ according as $j\equiv 0\,(\mod{\frac{q^m-1}{N}})$ or $j\not\equiv 0\,(\mod{\frac{q^m-1}{N}})$, continuing from \eqref{eq:hypfir2}, 
we have 
\[
S_{a,b}=\frac{1}{N}\sum_{j=0}^{N-1}
\sum_{\ell\in I}G_{q^m}(\chi_{N}^{-j})^2\chi_{N}^{j}(ab\omega^\ell). 
\]
On the other hand, by \eqref{eq:defIuse}, we have $
\sum_{i\in I}\chi_N^j(\omega^i)=\frac{G_{q^m}(\chi_N^j)}{\delta q}$. 
Hence, we have 
\begin{align*}
S_{a,b}-\frac{|I|}{N}
=\,&\frac{1}{\delta qN} \sum_{j=1}^{N-1}\chi_{N}^{j}(ab)G_{q^m}(\chi_{N}^{-j})G_{q^m}(\chi_{N}^{-j})G_{q^m}(\chi_{N}^{j})\nonumber\\
=\,& \frac{q^{m-1}}{\delta N} \sum_{j=1}^{N-1}\chi_{N}^{j}(\omega^{a})G_{q^m}(\chi_{N}^{-j})\\
=\,&\frac{q^{m}}{N}\sum_{j=1}^{N-1}\sum_{\ell \in I}\chi_N^{-j}(\omega^\ell)\chi_N^{j}(ab)
=-\frac{q^m|I|}{N}+\begin{cases}
q^m,& \text{ if }  \log_{\omega}(ab)\in I\, (\mod{N}), \\
0,& \text{ if } \log_{\omega}(ab)\not\in I\, (\mod{N}). 
\end{cases} 
\end{align*}
Thus $\psi_{a,b}(H)=S_{a,b}$, $(0,0)\neq (a,b)\in \F_{q^m}\times \F_{q^m}$, take exactly two distinct values $-\frac{(q^m-1)|I|}{N}$ and 
$-\frac{(q^m-1)|I|}{N}+q^m$. Therefore $\Cay(\F_{q^m}\times \F_{q^m}, H)$ is strongly regular. The parameters of $\Cay(\F_{q^m}\times \F_{q^m}, H)$ can be computed in a straightforward way. We omit the details. \qed
\begin{remark}{\em  
Under the assumptions of Proposition~\ref{theorem:main2}, set 
\[
H':=H\cup \{(0,x)\,|\,x \in \F_{q^m}^\ast\}\cup \{(x,0)\,|\,x \in \F_{q^m}^\ast\}. 
\]
Then, $\Cay(\F_{q^m}\times \F_{q^m}, H')$ is also strongly regular. Furthermore, if $\Cay(\F_{q^m}, C_0^{(N,q^m)})$ is a cyclotomic strongly regular graph in the subfield case
with $N=\frac{q^m-1}{q-1}$, we have \[
H'=\{(0,0)\neq (x,y)\in \F_{q^m}\times \F_{q^m}\,|\,\Tr_{q^m/q}(xy)=0\},
\] where $\Tr_{q^m/q}(xy)$ is a nondegenerate hyperbolic quadratic form from $\F_{q^m}\times \F_{q^m}$ to $\F_{q}$. Hence, it is appropriate to call the lifting construction in Proposition~\ref{theorem:main2} a {\it hyperbolic} type lifting construction.}
\end{remark}
We now apply Propositions~\ref{theorem:main1} and ~\ref{theorem:main2} to the known cyclotomic strongly regular graphs. 
We first apply the two propositions to the semi-primitive examples. In this case, 
we have $|I|=1$.     

\begin{corollary}\label{cor:semi}
Let $p$ be a prime, $N\ge 2$, $q^m=p^{2js}$, where  
$s\ge 2$, $N\,|\,(p^j+1)$, and $j$ is the smallest such positive integer. For $\epsilon\in \{1,-1\}$,  there exists an $(n^2,r(n-\epsilon),\epsilon n+r^2-3\epsilon r,r^2-\epsilon r)$ strongly regular Cayley graph with $n=q^m$ and $r=(q^m-1)/N$.
\end{corollary}
%
Next we apply Propositions~\ref{theorem:main1} and \ref{theorem:main2}  to the subfield examples. In this case, 
we have $N=\frac{q^{m}-1}{q^d-1}$ and $|I|=\frac{q^{m-d}-1}{q^d-1}$, where $d\,|\,m$. 
\begin{corollary}\label{cor:singer}
Let $q$ be a prime power and $m\ge 3$ a positive integer.  
For $\epsilon\in \{1,-1\}$,  there exists an $(n^2,r(n-\epsilon),\epsilon n+r^2-3\epsilon r,r^2-\epsilon r)$ strongly regular Cayley graph with 
$n=q^m$ and $r=q^{m-d}-1$.
\end{corollary}
When $d=1$, the strongly regular graphs obtained in Corollary~\ref{cor:singer} were already known~\cite{Ma}. 

Finally, we apply Propositions~\ref{theorem:main1} and \ref{theorem:main2} to the eleven sporadic  examples of cyclotomic strongly regular graphs. In this case,  the values of $k:=|I|$ are 
given in \cite[Table II]{SW}.  
\begin{corollary}\label{cor-spor}
For $\epsilon\in \{1,-1\}$, there exists an $(n^2,r(n-\epsilon),\epsilon n+r^2-3\epsilon r,r^2-\epsilon r)$ strongly regular Cayley graph with 
$n=q^m$ and $r=k(q^m-1)/N$ in each of the following cases:  
\begin{eqnarray*}
(q^m,N,k)&=&(3^5,11,5),(5^9,19,9),(3^{12},35,17),(7^9,37,9),(11^7,43,21),(17^{33},67,33)\\
& &(3^{53},107,53),(5^{18},133,33),(41^{81},163,81),(3^{144},323,161),(5^{249},499,249). 
\end{eqnarray*}
\end{corollary}

\section{Halving the connection sets $E$ and $H$ and their complements}

In a couple of recent papers \cite{FMX, BLMX}, motivated by existence questions concerning finite geometric objects such as $m$-ovoids and $i$-tight sets,  we used a certain partition of the Singer difference set to construct strongly regular Cayley graphs with special properties which give the desired $m$-ovoids and $i$-tight sets. We now realize that the constructions can be done in a more general setting, namely, we can do the construction by partitioning a subdifference set of the Singer difference set in a certain way. In the case where the cyclotomic strongly regular graph comes from a subfield, the subdifference set of the Singer difference set is actually a Singer difference set; so in this case, we recover the previous constructions. We will also use a certain partition of the complement of a subdifference of the Singer difference sets to construct more strongly regular Cayley graphs.

Assume that $N\geq 2$ is odd, $N\,|\,\frac{q^m-1}{q-1}$, and $\Cay(\F_{q^m},C_0^{(N,q^m)})$ is strongly regular. Let $I$ be the corresponding subdifference set in $\Z_N$ defined in (\ref{eq:subdi}). Let $S_1, S_2$ be a partition of $I$ and 
let  
$S_i'\equiv 2^{-1}S_i\,(\mod{N})$ and $S_i''\equiv 2^{-1}S_i'\,(\mod{N})$ for 
$i=1,2$. Define 
\begin{equation}\label{eq:defX1}
X:=2S_1'' \cup (2S_2''+N)\,(\mod{2N}). 
\end{equation}
Let $J_1:=\{0,3\}$ and 
$J_2:=\{1,2\}$. Define \begin{equation}\label{eq:defI}
Y:=\{Ni+4j \hspace{-0.2cm}\pmod{4N}: (i,j)\in (J_1\times S_1'') \cup (J_2\times S_2'')\}. 
\end{equation}
It is clear that $X\equiv 2^{-1}I\,(\mod{N})$ and $Y\equiv I\,(\mod{N})$. 

Similarly,  let $T_1, T_2$ be a partition of  $\Z_N\setminus I$ and let  
$T_i'\equiv 2^{-1}T_i\,(\mod{N})$ and $T_i''\equiv 2^{-1}T_i'\,(\mod{N})$ for 
$i=1,2$. 
Define 
\begin{equation}\label{eq:defX2}
\widehat{X}:=2T_1'' \cup (2T_2''+N)\,(\mod{2N}). 
\end{equation} 
Furthermore, define 
\begin{equation}\label{eq:defI2}
\widehat{Y}:=\{Ni+4j \hspace{-0.2cm}\pmod{4N}: (i,j)\in (J_1\times T_1'') \cup (J_2\times T_2'')\}, 
\end{equation}
where $J_1=\{0,3\}$ and $J_2=\{1,2\}$.

\subsection{Decompositions of  $\Cay(\F_{q^{2m}}, E)$ and its complement}\label{sec:minus}
In this subsection, we always assume that $q^m\equiv 3\,(\mod{4})$. We will consider decompositions of  $\Cay(\F_{q^{2m}}, E)$ and its complement, where $E$ is defined in (\ref{ellip}). We define
\begin{equation}\label{eq:defD-}
E_1:=\bigcup_{i \in Y}C_i^{(4N,q^{2m})}, 
\end{equation}
where $C_i^{(4N,q^{2m})}:=\gamma^i \langle \gamma^{4N}\rangle$, $\gamma$ is a primitive element of $\F_{q^{2m}}$, and $Y$ is defined in \eqref{eq:defI}. Since $Y\equiv I\pmod N$, we see that $E_1$ is a subset of $E$, and $|E_1|=|E|/2$. 
The (additive) character values of $E_1$ are given by the following lemma. 
\begin{lemma}\label{rem:quad-}
Let $\psi_{\F_{q^{2m}}}$ and $\psi_{\F_{q^m}}$ be the canonical additive characters of $\F_{q^{2m}}$ and $\F_{q^{m}}$, respectively. 
For $a\in \Z_{4N}$, define $b\equiv 4^{-1}a\,(\mod{N})$ and $c\equiv 2b\,(\mod{2N})$. Then, 
\begin{align}
\psi_{\F_{q^{2m}}}(\gamma^a E_1)=&\,\frac{\rho_p \delta_a q^m}{2G_{q^m}(\eta)}\left(2\psi_{\F_{q^m}}(\omega^c\bigcup_{t\in X}C_t^{(2N,q^m)})-\psi_{\F_{q^m}}
(\omega^c \bigcup_{t\in 2^{-1}I}C_t^{(N,q^m)})\right)\nonumber\\
&\, \hspace{3.0cm}+\frac{(q^m-1)|I|}{2N}-\left\{
\begin{array}{ll}
\frac{q^m}{2}, & \mbox{ if $c \in 2^{-1}I\,(\mod{N})$,}\\
0, & \mbox{ otherwise,}
 \end{array}
\right.\label{eq:DD-}
\end{align}
where $\delta_a=1$ or $-1$ depending on whether {$a\equiv 0,N\,(\mod{4})$ or $a\equiv 2,3N\,(\mod{4})$,} and $\rho_p=1$ or $-1$ depending on
whether $p\equiv 7\,(\mod{8})$ or $p\equiv 3\,(\mod{8})$. Furthermore, 
$\eta$ is the quadratic character of $\F_{q^m}$. 
\end{lemma}
This lemma is a common generalization of the results in \cite{BLMX} and \cite{Mo}. Its proof is the same as those in \cite{BLMX} and \cite{Mo}. We therefore omit the proof. 

Next, we consider a decomposition of the complement of $\Cay(\F_{q^{2m}}, E)$. Let 
\begin{equation}\label{eq:defD-2}
E_2:=\bigcup_{i \in \widehat{Y}}C_i^{(4N,q^{2m})}, 
\end{equation}
where $\widehat{Y}$ is defined in \eqref{eq:defI2}. 
The (additive) character values of $E_2$ are given by the following lemma. 
\begin{lemma}\label{rem:quad-2}
With the same notation as in Lemma~\ref{rem:quad-}, 
\begin{align}
\psi_{\F_{q^{2m}}}(\gamma^a E_2)=&\,\frac{\rho_p \delta_a q^m}{2G_{q^m}(\eta)}\left(2\psi_{\F_{q^m}}(\omega^c\bigcup_{t\in \widehat{X}}C_t^{(2N,q^m)})-\psi_{\F_{q^m}}
(\omega^c \bigcup_{t\in \Z_N\setminus 2^{-1}I}C_t^{(N,q^m)})\right)\nonumber\\
&\, \hspace{3.0cm}+\frac{(q^m-1)(N-|I|)}{2}-\left\{
\begin{array}{ll}
0, & \mbox{ if $c \in 2^{-1}I\,(\mod{N})$,}\\
\frac{q^m}{2}, & \mbox{ otherwise.}
 \end{array}
\right.\label{eq:DD-2}
\end{align}
\end{lemma}
\begin{remark}\label{rem:chara}{\em 
\begin{itemize}
\item[(i)] If  $X$ defined in \eqref{eq:defX1} satisfies that   
\begin{equation}\label{eq:3-chara}
2\psi_{\F_{q^m}}(\omega^c \bigcup_{t\in X}C_{t}^{(2N,q^m)})-
\psi_{\F_{q^m}}(\omega^c \bigcup_{t\in 2^{-1}I}C_{t}^{(N,q^m)})=\left\{
\begin{array}{ll}
\pm G_{q^m}(\eta), & \mbox{ if $c\in 2^{-1}I\, (\mod{N})$,}\\
0, & \mbox{ otherwise, }
 \end{array}
\right. 
\end{equation}
substituting \eqref{eq:3-chara} into \eqref{eq:DD-}, we find that the nontrivial additive character values of $E_1$ take two distinct values $\frac{(q^m-1)|I|}{2N}$ and $\frac{(q^m-1)|I|}{2N}-q^m$, implying that $\Cay(\F_{q^{2m}}, E_1)$ is strongly regular. 
\item[(ii)] If 
$\widehat{X}$ defined in \eqref{eq:defX2} satisfies that   
\begin{equation}\label{eq:3-chara2}
2\psi_{\F_{q^m}}(\omega^c \bigcup_{t\in \widehat{X}}C_{t}^{(2N,q^m)})-
\psi_{\F_{q^m}}(\omega^c \bigcup_{t\in \Z_N\setminus 2^{-1}I}C_{t}^{(N,q^m)})=\left\{
\begin{array}{ll}
0, & \mbox{ if $c\in 2^{-1}I\, (\mod{N})$,}\\
\pm G_{q^m}(\eta), & \mbox{ otherwise, }
 \end{array}
\right. 
\end{equation}
substituting \eqref{eq:3-chara2} into \eqref{eq:DD-2}, we find that the nontrivial additive character values of $E_2$ take two distinct values $\frac{(q^m-1)(N-|I|)}{2N}$ and $\frac{(q^m-1)(N-|I|)}{2N}-q^m$, implying that $\Cay(\F_{q^{2m}}, E_2)$ is strongly regular. 
\end{itemize}}
\end{remark}
\subsection{Decompositions of  $\Cay(\F_{q^{m}}\times \F_{q^{m}}, H)$ 
and its complement}
In this subsection, we assume that $q^m\equiv 1\,(\mod{4})$, $N$ is an odd divisor of $\frac{q^m -1}{q-1}$,  and $\gcd{(N,\frac{q^m-1}{N})}=1$. 
Define 
\begin{equation}\label{def:D+}
H_1:=\{(xy,xy^{-1}z\omega^{\ell})\,|\,x\in C_0^{(N,q^m)},y\in C_0^{(\frac{q^m-1}{N},q^m)},z\in C_0^{(4N,q^m)},\ell \in Y\}\subseteq
\F_{q^m}\times \F_{q^m}, 
\end{equation}
where $\omega$ is a primitive element of $\F_{q^m}$ defined in Subsection~\ref{sec:two} and $Y$ is defined in \eqref{eq:defI}.  In the definition of $H_1$, since $x^2z\omega^\ell \in \bigcup_{\ell \in I}C_\ell^{(N,q^m)}$, we see that $H_1$ is a subset of $H$. Moreover, $|H_1|=|H|/2.$ The (additive) character values of $H_1$ are given in the following lemma. 
\begin{lemma}\label{mainconstruction2+} 
Let $\psi_{a,b}$ be an additive character of $\F_{q^m}\times \F_{q^m}$ defined in \eqref{fieldchara22} and $\psi_{\F_{q^m}}$ be the canonical additive character of $\F_{q^m}$. 
For $(a,b)\in \F_{q^m}\times \F_{q^m}\setminus \{(0,0)\}$ with $ab=0$, it holds that $\psi_{a,b}(H_1)
=-(q^m-1)|I|/2N$. For $(a,b)\in \F_{q^m}\times \F_{q^m}\setminus \{(0,0)\}$ with $ab\not=0$, it holds that   
\begin{align}
\psi_{a,b}(H_1)=&\,\frac{\eta(2\omega^c)G_{q^m}(\eta)\delta_{a,b} 
}{2}
\left(2\psi_{\F_{q^m}}(\omega^c\bigcup_{t\in X}C_t^{(2N,q^3)})-\psi_{\F_{q^m}}
(\omega^c \bigcup_{t\in 2^{-1}I}C_t^{(N,q^m)})\right)\nonumber\\
&\, \hspace{3.0cm}-\frac{(q^m-1)|I|}{2N}+\left\{
\begin{array}{ll}
\frac{q^m}{2}, & \mbox{ if $c \in 2^{-1}I\,(\mod{N})$,}\\
0, & \mbox{ otherwise,}
 \end{array}
\right.\label{eq:DD+}
\end{align}
where $c$ is defined by  $\omega^c=(ab)^{\frac{N+1}{2}}$ and $\delta_{a,b}=1$ or $-1$ depending on whether {$\log_{\omega}(a^{-1}b)\equiv 0,N\,(\mod{4})$ or $\log_{\omega}(a^{-1}b)\equiv 2,3N\,(\mod{4})$. } 
Furthermore, 
$\eta$ is the quadratic character of $\F_{q^m}$. 
\end{lemma}
This lemma is a generalization of \cite[Theorem~4.1]{FMX}. Since the proof is similar to that of \cite[Theorem~4.2]{FMX}, we omit the proof. 

Next, we consider a decomposition of the complement of $\Cay(\F_{q^{m}}\times \F_{q^m}, H)$. Define 
\begin{equation}\label{def:D+}
H_2:=\{(xy,xy^{-1}z\omega^{\ell})\,|\,x\in C_0^{(N,q^m)},y\in C_0^{(\frac{q^m-1}{N},q^m)},z\in C_0^{(4N,q^m)},\ell \in \widehat{Y}\}\subseteq
\F_{q^m}\times \F_{q^m}, 
\end{equation}
where $\widehat{Y}$ is defined in \eqref{eq:defI2}. 
The character values of $H_2$ are given in the following lemma. 
\begin{lemma}\label{mainconstruction3+} 
For $(a,b)\in \F_{q^m}\times \F_{q^m}\setminus \{(0,0)\}$ with $ab=0$, it holds that $\psi_{a,b}(H_2)
=-(q^m-1)(N-|I|)/2N$. For  $(a,b)\in \F_{q^m}\times \F_{q^m}\setminus \{(0,0)\}$ with $ab\not=0$,  it holds that   
\begin{align}
\psi_{a,b}(H_2)=&\,\frac{\eta(2\omega^c)G_{q^m}(\eta)\delta_{a,b} 
}{2}
\left(2\psi_{\F_{q^m}}(\omega^c\bigcup_{t\in \widehat{X}}C_t^{(2N,q^3)})-\psi_{\F_{q^m}}
(\omega^c \bigcup_{t\in \Z_N\setminus 2^{-1}I}C_t^{(N,q^m)})\right)\nonumber\\
&\, \hspace{3.0cm}-\frac{(q^m-1)(N-|I|)}{2N}+\left\{
\begin{array}{ll}
0, & \mbox{ if $c \in 2^{-1}I\,(\mod{N})$,}\\
\frac{q^m}{2}, & \mbox{ otherwise, }
 \end{array}
\right.\label{eq:DD+2}
\end{align}
where $c$ is defined by $\omega^c=(ab)^{\frac{N+1}{2}}$. 
\end{lemma}
\begin{remark}\label{re:cha:hyp}{\em 
Similarly to Remark~\ref{rem:chara}, 
if  the set $X$ defined in \eqref{eq:defX1} satisfies \eqref{eq:3-chara}, 
the nontrivial additive character values of $H_1$ take two distinct values $-\frac{(q^m-1)|I|}{2N}$ and $-\frac{(q^m-1)|I|}{2N}+q^m$, implying that $\Cay(\F_{q^m}\times \F_{q^m}, H_1)$ is strongly regular. 
Also, if 
$\widehat{X}$ defined in \eqref{eq:defX2} satisfies 
\eqref{eq:3-chara2}, 
then the nontrivial additive character values of $H_2$ take two distinct values $-\frac{(q^m-1)(N-|I|)}{2N}$ and $-\frac{(q^m-1)(N-|I|)}{2N}+q^m$, implying that $\Cay(\F_{q^m}\times \F_{q^m}, H_2)$ is strongly regular.} 
\end{remark}

\section{Partition of subdifference sets and their complements in semi-primitive case}
In this section, we consider a partition of the subdifference sets $I$ in semi-primitive case. We will use the same notation as in Section 4.
We assume that $N$ is odd and $q^{m}=p^{2js}$, where $p$ is a prime, 
$s\ge 2$, $N\,|\,(p^j+1)$, and $j$ is the smallest 
such positive integer. In this case, $\Cay(\F_{q^{m}},C_0^{(N,q^m)})$ is strongly regular and we have $I=\{0\}$~\cite{SW}. 
Furthermore, by Theorem~\ref{thm:semiprim}, the Gauss sums with respect to multiplicative characters of exponent $N$ of $\F_{q^m}$ can be explicitly evaluated as 
\begin{equation}\label{eq:Gsemi}
G_{q^m}(\chi_N^i)=(-1)^{s-1}\sqrt{q^m}, \quad 1\le i\le N-1. 
\end{equation}
\begin{theorem}\label{thm:semi-}
With the same notation as in Section 4, under the above assumptions, the partition $(S_1,S_2)=(\{0\},\emptyset)$ of $I$ satisfies the 
condition~\eqref{eq:3-chara} of Remark~\ref{rem:chara} (i). 
\end{theorem}
\proof 
By the definition~\eqref{eq:defX1} of $X$, we have $X=\{0\}$.  
Write
\[
P_c:=2\psi_{\F_{q^m}}(\omega^c \bigcup_{t\in X}C_t^{(2N,q^m)})-\psi_{\F_{q^m}}(\omega^c \bigcup_{t\in 2^{-1}I}C_t^{(N,q^m)}). 
\]
By \eqref{eq:3-chara}, we need to prove that \[
P_c=\left\{
\begin{array}{ll}
\pm G_{q^m}(\eta), & \mbox{ if $c\equiv 0\, (\mod{N})$,}\\
0, & \mbox{ otherwise,}
 \end{array}
\right. 
\]
where $\eta$ is the quadratic  character of $\F_{q^m}$. 
Let $\chi_{N}$ be a multiplicative character of {order $N$} of $\F_{q^m}$. 
By the orthogonality of characters, we have
\begin{align}
P_c=&\,\frac{1}{N}\left(\sum_{i=0}^{N-1}\sum_{j=0,1}\sum_{t\in X}G_{q^m}(\chi_{N}^i\eta^j)\chi_{N}^{-i}\eta^{j}(\omega^{c+t})-\sum_{i=0}^{N-1}\sum_{t\in 2^{-1}I}G_{q^m}(\chi_{N}^i)\chi_{N}^{-i}(\omega^{c+t})\right) \nonumber\\
=&\,\frac{1}{N}\sum_{i=0}^{N-1}\sum_{t\in X}G_{q^m}(\chi_{N}^i\eta)\chi_{N}^{-i}\eta(\omega^{c+t}). \label{eq:hiki}
\end{align}
By the Davenport-Hasse product formula and \eqref{eq:Gsemi}, we have 
\[
G_{q^m}(\chi_{N}^i\eta)=\frac{G_{q^m}(\chi_N^{2i})G_{q^m}(\eta)}{G_{q^m}(\chi_{N}^i)}=G_{q^m}(\eta). 
\]
On the other hand, by the definition of $X$, we have 
$\sum_{t\in X}\chi_{N}^{-i}\eta(\omega^{t})=1$. 
Continuing from \eqref{eq:hiki}, we have 
\[
P_c=\frac{\eta(\omega^{c})G_{q^m}(\eta)}{N}\sum_{i=0}^{N-1}\chi_{N}^{-i}(\omega^{c})
=
\left\{
\begin{array}{ll}
\eta(\omega^{c})G_{q^m}(\eta), & \mbox{ if $c \equiv 0\,(\mod{N})$,}\\
0, & \mbox{ otherwise.}
 \end{array}
\right. \nonumber
\]
This completes the proof of the theorem. 
\qed

\vspace{0.3cm}
Similarly to the theorem above, we have the following. 
\begin{theorem}\label{thm:spoc2}
With the notations above, 
the partition $(T_1, T_2)=(\Z_N\setminus \{0\},\emptyset)$ of $\Z_N\setminus I$ satisfies the 
condition~\eqref{eq:3-chara2}  of Remark~\ref{rem:chara} (ii). 
\end{theorem}
Since $q^m=p^{2js}$, we have $q^m\equiv 1\,(\mod{4})$.  We can only apply the lifting construction of hyperbolic type. 
By Lemma~\ref{mainconstruction2+}, Remark~\ref{re:cha:hyp} and Theorem~\ref{thm:semi-},  we obtain the following. 
\begin{corollary}\label{cor-semi-}
Let $N$ be odd and $q^{m}=p^{2js}$, where $p$ is a prime, 
$s\ge 2$, $N\,|\,(p^j+1)$, and $j$ is the smallest 
such positive integer. Assume that $\gcd{(N,\frac{q^m-1}{N})}=1$.  Then, 
there exists a $(q^{2m},r(q^{m}-1),q^m+r^2-3 r,r^2- r)$ strongly regular Cayley graph, where $r=(q^m-1)/2N$.   
\end{corollary}

\vspace{0.3cm}
Similarly to the corollary above, by Lemma~\ref{mainconstruction3+}, Remark~\ref{re:cha:hyp} and Theorem~\ref{thm:spoc2}, we obtain the following corollary. 
\begin{corollary}\label{cor-semi-2}
Let $N$ be odd and $q^{m}=p^{2js}$, where $p$ is a prime, 
$s\ge 2$, $N\,|\,(p^j+1)$, and $j$ is the smallest 
such positive integer. Assume that $\gcd{(N,\frac{q^m-1}{N})}=1$. Then, there exists a $(q^{2m},r(q^{m}-1),q^m+r^2-3 r,r^2- r)$ strongly regular Cayley graph, where  $r=(N-1)(q^m-1)/2N$. 
\end{corollary}

\section{Partition of subdifference sets and their complements in sporadic case}
In this section, we consider partitions of the subdifference set $I$ and its complement in the sporadic case. 
\begin{theorem}\label{thm:spo}
Assume that $N\geq 2$ is odd, $N\;|\;\frac{q^m-1}{q-1}$, $\Cay(\F_{q^m},C_0^{(N,q^m)})$ is strongly regular and $-2\in \langle p\rangle\,(\mod{N})$, where $p$ is the characteristic of $\F_{q^m}$. Let $I$ be the corresponding subdifference set defined in (\ref{eq:subdi}). Then, the partition $(S_1,S_2)=(I,\emptyset)$ of $I$ satisfies the 
condition~\eqref{eq:3-chara} of Remark~\ref{rem:chara} (i). 
\end{theorem}
\proof
Let $I'\equiv 4^{-1}I\,(\mod{N})$. 
By the definition~\eqref{eq:defX1} of $X$, we have $X\equiv 2I'\,(\mod{2N})$.  
Write
\[
P_c:=2\psi_{\F_{q^m}}(\omega^c \bigcup_{t\in X}C_t^{(2N,q^m)})-\psi_{\F_{q^m}}(\omega^c \bigcup_{t\in 2^{-1}I}C_t^{(N,q^m)}). 
\]
Let $\chi_{N}$ be a multiplicative character of $\F_{q^m}$ of order $N$ and 
$\eta$ be the quadratic  character of $\F_{q^m}$. Similarly to the proof of 
Theorem~\ref{thm:semi-}, we have 
\begin{align}
P_c=\frac{1}{N}\sum_{i=1}^{N-1}\sum_{t\in X}G_{q^m}(\chi_{N}^i\eta)\chi_{N}^{-i}\eta(\omega^{c+t})+\frac{\eta(\omega^c)G_{q^m}(\eta)|I|}{N}. \label{eq:hiki2}
\end{align}
By the Davenport-Hasse product formula, we have 
\[
G_{q^m}(\chi_{N}^i\eta)=\frac{G_{q^m}(\chi_N^{2i})G_{q^m}(\eta)}{G_{q^m}(\chi_{N}^i)}. 
\]
On the other hand, by  \eqref{eq:sum},  for $i\not=0$ 
\[
\sum_{t\in X}\chi_{N}^{-i}\eta(\omega^{t})=\sum_{t\in 2I'}\chi_{N}^{-i}(\omega^{t})=\frac{G_{q^m}(\chi_N^{-2^{-1}i})}{\delta q}. 
\]
Substituting these into \eqref{eq:hiki2}, we have 
\begin{equation}\label{eq:subdi2}
P_c=\frac{\eta(\omega^{c})G_{q^m}(\eta)}{\delta q N}\sum_{i=1}^{N}\frac{G_{q^m}(\chi_N^{2i})G_{q^m}(\chi_N^{-2^{-1}i})}{G_{q^m}(\chi_{N}^i)}\chi_{N}^{-i}(\omega^{c})+\frac{\eta(\omega^c)G_{q^m}(\eta)|I|}{N}. 
\end{equation}
By the assumption that $-2\in \langle p\rangle \,(\mod{N})$, we have 
$G_{q^m}(\chi_N^{-2^{-1}i})=G_{q^m}(\chi_N^i)$. Therefore, continuing from 
\eqref{eq:subdi2}, we have 
\begin{align*}
P_c=&\,\frac{\eta(\omega^{c})G_{q^m}(\eta)}{\delta q N}\sum_{i=1}^{N}G_{q^m}(\chi_N^{2i})\chi_{N}^{-i}(\omega^{c})+\frac{\eta(\omega^c)G_{q^m}(\eta)|I|}{N}\\
=&\,\frac{\eta(\omega^{c})G_{q^m}(\eta)}{N}\left(\sum_{i=1}^{N}\sum_{t\in I}\chi_N^{2i}(\omega^t)\chi_{N}^{-i}(\omega^{c})+|I|\right)\\
=&\,\frac{\eta(\omega^{c})G_{q^m}(\eta)}{N}\sum_{i=0}^{N}\sum_{t\in I}\chi_N^{2i}(\omega^t)\chi_{N}^{-i}(\omega^{c})=
\left\{
\begin{array}{ll}
\eta(\omega^{c})G_{q^m}(\eta), & \mbox{ if $c \in 2I\,(\mod{N})$,}\\
0, & \mbox{ otherwise.}
 \end{array}
\right.
\end{align*}
Since the subdifference set $I$ is invariant under the multiplication by $p$ modulo $N$, by the assumption that $-2^{-1}\in \langle p\rangle \,(\mod{N})$, the condition 
$c\in 2I\,(\mod{N})$ is equivalent to that $c\in 2^{-1}I\,(\mod{N})$.  
This completes the proof of the theorem. 
\qed

\vspace{0.3cm}
Similarly to the theorem above, we have the following. 
\begin{theorem}\label{thm:spoc}
Assume that $N\geq 2$ is odd, $N\;|\;\frac{q^m-1}{q-1}$,  $\Cay(\F_{q^m},C_0^{(N,q^m)})$ is strongly regular and $-2\in \langle p\rangle\,(\mod{N})$. 
Then the partition $(T_1, T_2)=(\Z_N\setminus I,\emptyset)$ of $\Z_N\setminus I$ satisfies the 
condition~\eqref{eq:3-chara2}  of Remark~\ref{rem:chara} (ii). 
\end{theorem}
There are ten sporadic examples of cyclotomic strongly regular graphs  satisfying the condition $-2\in \langle p\rangle\,(\mod{N})$. In particular, when 
$q^m\equiv 3\,(\mod{4})$, we obtain the following result. 
\begin{corollary}\label{cor-spor}
There exists a $(q^{2m},r(q^{m}+1),-q^m+r^2+3 r,r^2+ r)$ strongly regular Cayley graph with $r=k(q^m-1)/2N$  in each of the following cases:  
\[
(q^m,N,k)=(3^5,11,5),(11^7,43,21),(3^{53},107,53). 
\]
\end{corollary}
\proof 
It is clear that $-2\in \langle p\rangle\,(\mod{N})$ in these cases. Then, by applying Lemma~\ref{rem:quad-}, Remark~\ref{rem:chara}~(i)  and  Theorem~\ref{thm:spo} to these examples, the corollary now follows. 
\qed

\vspace{0.3cm}
Similarly to the corollary above, by applying Lemma~\ref{rem:quad-2},  Remark~\ref{rem:chara}~(ii) and Theorem~\ref{thm:spoc}  to the three sporadic cyclotomic strongly regular graphs in Corollary~\ref{cor-spor}, we obtain the following. 
\begin{corollary}\label{cor-spor2}
There exists a $(q^{2m},r(q^{m}+1),-q^m+r^2+3 r,r^2+ r)$ strongly regular Cayley graph with $r=(N-k)(q^m-1)/2N$  in each of the following cases:  
\[
(q^m,N,k)=(3^5,11,5),(11^7,43,21),(3^{53},107,53). 
\]
\end{corollary}

\vspace{0.3cm}
In the case where $(q^m,N,k)=(7^9,37,9)$, 
the condition that $-2\in \langle p\rangle\,(\mod{N})$ is not satisfied. 
We checked by computer that there is no partition of the subdifference set $I$ satisfying the condition~\eqref{eq:3-chara}. On the other hand, we checked that there is a  
partition of $\Z_{37}\setminus I$ satisfying the condition~\eqref{eq:3-chara2}: $T_1= 2I$ and  
$T_2=\Z_{37}\setminus (I\cup 2I)$. 
Hence, 
we have the following corollary. 
\begin{corollary}\label{cor-spor}
There exists a $(q^{2m},r(q^m+1),-q^m+r^2+3 r,r^2+ r)$ strongly regular Cayley graph with $r=(N-k)(q^m-1)/2N$ in the case where $(q^m,N,k)=(7^9,37,9)$. 
\end{corollary}
Next, we consider the case where $q^m\equiv 1\,(\mod{4})$. 
\begin{corollary}\label{cor-spor3}
There exists a $(q^{2m},r(q^m-1),q^m+r^2-3r,r^2- r)$ strongly regular Cayley graph with $r=k(q^m-1)/2N$ in each of the following cases:  
\begin{align*}
(q^m,N,k)=&\,(3^{12},35,17),(5^9,19,9),(17^{33},67,33), (5^{18},133,33),\\
&\,(41^{81},163,81),(3^{144},323,161),(5^{249},499,249). 
\end{align*}
\end{corollary}
\proof 
It is clear that  $\gcd{(N,\frac{q^m-1}{N})}=1$ and $-2\in \langle p\rangle\,(\mod{N})$ in these cases. 
Then, by applying Lemma~\ref{mainconstruction2+}, Remark~\ref{re:cha:hyp} and  Theorem~\ref{thm:spo} to these examples, the corollary now follows. 
\qed

\vspace{0.3cm}
Similarly to the corollary above, by applying Lemma~\ref{mainconstruction3+}, Remark~\ref{re:cha:hyp} and Theorem~\ref{thm:spoc} to these examples, we obtain the following corollary. 
\begin{corollary}\label{cor-spor4}
There exists a $(q^{2m},r(q^{m}-1),q^m+r^2-3 r,r^2- r)$ strongly regular Cayley graph with $r=(N-k)(q^m-1)/2N$  in each of the following cases:  
\begin{align*}
(q^m,N,k)=&\,(3^{12},35,17),(5^9,19,9),(17^{33},67,33), (5^{18},133,33),\\
&\,(41^{81},163,81),(3^{144},323,161),(5^{249},499,249). 
\end{align*}
\end{corollary}
\section{Partitions of subdifference sets and their complements in the subfield case}
In this section, we consider partitions of the subdifference set $I$ and its complement in subfield case. We assume that $m$ is odd and $N=\frac{q^m-1}{q-1}$. 
In this case, $\Cay(\F_{q^{m}},C_0^{(N,q^m)})$ is strongly regular and we have
\begin{equation}\label{eq:sub1}
I:=\{i\,(\mod{N}):\Tr_{q^m/q}(w^i)=0\}. 
\end{equation}
\subsection{A partition of the Singer difference set $I$ defined in (\ref{eq:sub1}) when $m=3$}
In the case where $m=3$, a partition of the Singer difference set $I$ satisfying the condition~\eqref{eq:3-chara} of Remark~\ref{rem:chara}~(i)  was found in \cite[Theorem~3.7]{FMX}. 
Regarding $\F_{q^3}$ as a $3$-dimensional vector space over $\F_q$, we  use $\F_{q^3}$ as the underlying vector space of $\PG(2,q)$. The points of $\PG(2,q)$ are $\langle{\omega^i}\rangle$, 
$0\le i\le N-1$, and the lines of $\PG(2,q)$ are 
\begin{equation}\label{eqn_Lu}
L_i:=\{\langle{x}\rangle:\,\Tr_{q^3/q}(\omega^i x)=0\},\, \, \, \, 0\le i\le N-1. 
\end{equation}
The Singer difference  set $I$ corresponds to the typical line $L_0$.  

Consider a nondegenerate quadratic form $f: \F_{q^3}\rightarrow \F_q$ defined by $f(x)=
\tr_{q^3/q}(x^2)$, which  
defines a conic $\cQ$ in $\PG(2,q)$ containing $q+1$ points. 
Consequently, each line $L$ of $\PG(2,q)$ meets $\cQ$ in $0$, $1$ or $2$ points. 
Consider the following subset of $\Z_N$: 
\begin{equation}\label{eqn_IQ}
I_\cQ:=\{i\, (\mod{N}):f(\omega^{i})=0\}=\{d_0,d_1,\ldots, d_{q}\},
\end{equation}
where the elements are numbered in any unspecified order. 
Thus,  $\cQ=\{\langle{\omega^{d_i}}\rangle:\,0\le i\le q\}$. Furthermore,  
by the definition of $f$ and $I$, $I_\cQ\equiv 2^{-1}I\pmod{N}$. 

For $d_0\in I_\cQ$, define
\[
{\mathcal X}:=\{\omega^{d_i}\Tr_{q^3/q}(\omega^{d_0+d_i}):\,1\le i\le q\}\cup\{2  \omega^{d_0}\}
\]
and 
\begin{equation}\label{eqn_defX2}
X:=\{\log_{\omega}(x)\, (\mod{2N}):\, x\in {\mathcal X}\}\subset \Z_{2N}. 
\end{equation}
Clearly, $|X|=|I_{\mathcal Q}|$ and $X\equiv I_\cQ\pmod{N}$. 
If we use any other $d_i$ instead of $d_0$ in the definition of ${\mathcal X}$, then the resulting set $X'$ satisfies that 
$X'\equiv X\, (\mod{2N})$ or $X'\equiv X+N\, (\mod{2N})$~\cite[Lemma~3.4]{FMX}. 


The set $X$ can be expressed as 
\begin{equation}\label{eqn_defX}
X=2S_1''\cup (2S_2''+N)\pmod{2N}
\end{equation}
for some $S_1'',S_2''\subseteq \Z_N$ with $|S_1''|+|S_2''|=q+1$. 
Define $S_i'\equiv 2S_i''\pmod{N}$ and $S_i\equiv 2S_i'\pmod{N}$ for 
$i=1,2$. 
Then,  
$S_1'\cup S_2'\equiv I_\cQ\pmod{N}$ and $S_1\cup S_2\equiv I\pmod{N}$, i.e., $X$ induces partitions of $I_\cQ$ and  $I$, respectively.  

\begin{theorem}\label{thm:sub6}{\em \cite[Theorem 3.7]{FMX}}
The set $X$ defined in \eqref{eqn_defX2} satisfies the condition \eqref{eq:3-chara}  of Remark~\ref{rem:chara} (i). 
\end{theorem}
As corollaries, we have the following. 
\begin{corollary}\label{cor:el}
For a prime power $q\equiv 3\,(\mod{4})$, there exists 
a $(q^{6},r(q^{3}+1),-q^3+r^2+3 r,r^2+ r)$ strongly regular Cayley graph, where $r=(q^2-1)/2$. 
\end{corollary}
\proof
By Lemma~\ref{rem:quad-}, Remark~\ref{rem:chara}~(i)  and  Theorem~\ref{thm:sub6}, the corollary now follows. 
\qed

\vspace{0.3cm} The connection set $E_1\subseteq \F_{q^6}$ of the strongly regular Cayley graph $\Cay(\F_{q^6}, E_1)$ obtained in Corollary~\ref{cor:el} 
corresponds to a $\frac{(q+1)}{2}$-ovoid in an elliptic quadric ${\mathcal Q}^-(5,q)$. See \cite{BLMX}. 
\begin{corollary}\label{cor:hyp}
For a prime power $q\equiv 5,9\,(\mod{12})$, there exists a  $(q^6,r(q^3-1),q^3+r^2-3r,r^2-r)$ strongly regular Cayley graph, where  $r=(q^2-1)/2$. 
\end{corollary}
\proof
It is clear that  $\gcd{(N,\frac{q^3-1}{N})}=1$ if $N=q^2+q+1$ and $q\equiv 5,9\,(\mod{12})$.  
Then, by Lemma~\ref{mainconstruction2+}, Remark~\ref{rem:chara}~(ii)  and  Theorem~\ref{thm:sub6}, the corollary now follows. 
\qed

\vspace{0.3cm} 
The connection set {$H_1\subseteq \F_{q^3}\times \F_{q^3}$} of the strongly regular Cayley graph {$\Cay(\F_{q^3}\times \F_{q^3}, H_1)$} obtained in Corollary~\ref{cor:hyp} 
corresponds to a $\frac{(q^2-1)}{2}$-tight set in a hyperbolic quadric ${\mathcal Q}^+(5,q)$. See \cite{DDMR,FMX}. 

It would be interesting to find a desired partition of $I$ when $m$ is odd and $m>3$. We leave this as an open problem.

\subsection{A partition of the complement of the Singer difference set with odd $m$}
In this section, we consider a partition of the complement of the Singer difference set $I\,(\mod{\frac{q^m-1}{q-1}})$, where  $m>1$ is an arbitrary odd integer. Note that the set $2^{-1}I\,(\mod{\frac{q^m-1}{q-1}})$ corresponds to a nondegenerate parabolic quadric ${\mathcal Q}(m-1,q)$ of $\PG(m-1,q)$. 

Let $N=\frac{q^m-1}{q-1}$ and $\omega$ be a primitive element of $\F_{q^m}$, where $q$ is an odd prime power and $m>1$ is an odd integer. 
Define 
\begin{align*}
A=&\,\{x\in \F_{q^m}^\ast\,|\,\Tr_{q^m/q}(x^{2})= 0\},\\
A_0=&\,\{x\in \F_{q^m}^\ast\,|\,\Tr_{q^m/q}(x^{2})\in C_0^{(2,q)}\},\\
A_1=&\,\{x\in \F_{q^m}^\ast\,|\,\Tr_{q^m/q}(x^{2})\in C_1^{(2,q)}\}.
\end{align*}
Let $a_1\in A$, and define $H_1=\{x\in \F_{q^m}^\ast\mid \Tr_{q^m/q}(a_1x)=0\}$. Note that $A$ represents a nondegenerate parabolic quadric of $\PG(m-1,q)$ and $H_1$ is a tangent hyperplane\footnote{Strictly speaking, we should say that {$H_1\cup\{0\}$ is a hyperplane.}} to $A$ at point $\langle a_1\rangle$. Thus $A\cap H_1$ is a cone of order one with vertex $\langle a_1\rangle$, and $|A\cap H_1|=q^{m-2}-1$. If $m=3$, we stop this process. Otherwise, we continue by choosing $a_2\in A\cap H_1$ such that $a_1,a_2$ are linearly independent over $\F_q$, and define $H_2=\{x\in \F_{q^m}^\ast\mid \Tr_{q^m/q}(a_2x)=0\}$. Then $H_1\cap H_2$ is a hyperplane of $H_1$. Note that $A\cap H_1$ represents a degenerate quadric (a cone of order one) in $H_1$, and $H_1\cap H_2$ contains the vertex $\langle a_1\rangle$, we see that $A\cap H_1\cap H_2$ is a cone of order two (cf. \cite{games}), and $|A\cap H_1\cap H_2|=q^{m-3}-1$. More generally, we define
\begin{align*}
H_{\ell}=&\,\{x\in \F_{q^m}^\ast\,|\,\Tr_{q^m/q}(xa_\ell)= 0\}, \, \quad a_\ell\in A \cap H_{1}\cap \cdots \cap H_{\ell-1},\\
H_{\ell,0}=&\,\{x\in \F_{q^m}^\ast\,|\,\Tr_{q^m/q}(xa_\ell)\in C_0^{(2,q)}\}, \, \quad a_\ell\in A \cap H_{1}\cap \cdots \cap H_{\ell-1},\\
H_{\ell,1}=&\,\{x\in \F_{q^m}^\ast\,|\,\Tr_{q^m/q}(xa_\ell)\in C_1^{(2,q)}\}, \, \quad a_\ell\in A \cap H_{1}\cap \cdots \cap H_{\ell-1},
\end{align*}
where $2\le \ell \le \frac{m-1}{2}$. We can always choose $a_1,\ldots, a_{\frac{m-1}{2}}$ so that they are linearly independent over $\F_q$. The reason is as follows: assume that $a_1,\ldots,a_{\ell-1}$ with $2\le \ell \le \frac{m-1}{2}$ are independent; since $a_1,\ldots,a_{\ell-1} \in A\cap H_{1}\cap \cdots \cap H_{\ell-1}$ and 
\begin{equation}\label{eq:numnu}
|A\cap H_{1}\cap \cdots \cap H_{\ell-1}|=q^{m-\ell}-1, 
\end{equation}
there are at least $m-\ell$ independent elements in $A\cap H_{1}\cap \cdots \cap H_{\ell-1}$ including $a_1,\ldots,a_{\ell-1}$; hence, we can choose an element $a_{\ell} \in A\cap H_{1}\cap \cdots \cap H_{\ell-1}$ so that $a_1,\ldots,a_\ell$ are independent over $\F_q$ whenever $\ell \le \frac{m-1}{2}$. 

Let $b$ be a fixed element of $(H_1\cap \cdots \cap H_{\frac{m-1}{2}})\setminus A$. 
Since $H_1\cap \cdots \cap H_{\frac{m-1}{2}}$ and  $A\cap H_1\cap \cdots \cap H_{\frac{m-1}{2}}$ correspond to a $\frac{(m-1)}{2}$-flat and 
 a $\frac{(m-3)}{2}$-flat, respectively, in $\PG(m-1,q)$, the set $(H_1\cap \cdots \cap H_{\frac{m-1}{2}})\setminus A$ can be represented as 
\[
(H_1\cap \cdots \cap H_{\frac{m-1}{2}})\setminus A=\{a_1x_1+\cdots+
a_\frac{m-1}{2} x_\frac{m-1}{2}+by\,|\,x_1,\ldots,x_{\frac{m-1}{2}}\in \F_q, y \in \F_{q}^\ast\}.
\]
Let $T_1=(A_0\cap H_{1,0})\cup (A_1\cap H_{1,1})$ and more generally
\[
T_\ell:=(A_0 \cap H_{1} \cap \cdots \cap H_{\ell-1} \cap H_{\ell,0})
\cup (A_1 \cap H_{1} \cap \cdots \cap H_{\ell-1} \cap H_{\ell,1}), 
\quad \, 2\le \ell \le \frac{m-1}{2}, 
\]
and 
\[
B:=\{a_1x_1+\cdots+
a_\frac{m-1}{2} x_\frac{m-1}{2}+by\,|\,x_1,\ldots,x_{\frac{m-1}{2}}\in \F_q, y \in C_0^{(2,q)}\}.
\]
Finally,  
define \[
D:=\left(\bigcup_{\ell=1}^{(m-1)/2}T_\ell\right)\cup B. 
\]
It is clear that 
\[
\omega^{N} T_\ell=(A_0 \cap H_{1} \cap \cdots \cap H_{\ell-1} \cap H_{\ell,1})
\cup (A_1 \cap H_{1} \cap \cdots \cap H_{\ell-1} \cap H_{\ell,0}) 
\] 
and \[
\omega^NB=\{a_1x_1+\cdots+
a_\frac{m-1}{2} x_\frac{m-1}{2}+by\,|\,x_1,\ldots,x_{\frac{m-1}{2}}\in \F_q, y \in C_1^{(2,q)}\}.
\]
Hence, $D\cap \omega^{N}D=\emptyset$ and 
$D\cup \omega^{N}D=\F_{q^m}^\ast \setminus A$. Thus, 
there exists a subset $\widehat{X}\subseteq \Z_{2N}$ such that 
$D=\bigcup_{t\in \widehat{X}}C_t^{(2N,q^m)}$ and $\widehat{X}\equiv \Z_N\setminus 2^{-1}I \,(\mod{N})$. The set $\widehat{X}$ induces a partition of  the complement of $2^{-1}I\,(\mod{N})$. 
\begin{theorem}\label{thm:coSing}
The set $\widehat{X}$ defined above satisfies the 
condition~\eqref{eq:3-chara2}  of Remark~\ref{rem:chara} (ii). 
\end{theorem}
To prove this theorem, we need the following lemmas.
\begin{lemma}\label{lem:coSing}
It holds that 
\[
\psi_{\F_{q^m}}(\omega^a 
T_\ell)= 
\left\{
\begin{array}{ll}
{\frac{(-1)^{i+\epsilon \frac{m-1}{2}}q^{\frac{m-1}{2}}(-1+(-1)^{j+\tau} G_q(\eta'))}{2},} & \mbox{ if $\omega^a \in A_i\cap H_1\cap \cdots \cap H_{\ell-1}\cap H_{\ell,j}$, $i,j=0,1,$ }\\
-\frac{q^{m-\ell-1}(q-1)}{2}, & \mbox{ if $\omega^a \in \langle a_1,\ldots,a_\ell\rangle\setminus \langle a_1,\ldots,a_{\ell-1}\rangle$, }\\
\frac{q^{m-\ell-1}(q-1)^2}{2}, & \mbox{ if $\omega^a \in \langle a_1,\ldots,a_{\ell-1}\rangle$, }\\
0, & \mbox{ otherwise, }
 \end{array}
\right.
\]
where $\eta'$ is the quadratic character of $\F_q$ and $\epsilon=0$ or $1$ according as $q\equiv 1\,(\mod{4})$ or $q\equiv 3\,(\mod{4})$. 
{Furthermore, $\tau$ is defined by $2\in C_\tau^{(2,q)}$.}  
\end{lemma}
The proof of this lemma is complicated. Therefore, we postpone the proof to the Appendix. 
\begin{lemma}\label{lem:coSing2}
It holds that  
\begin{align*}
\psi_{\F_{q^m}}(\omega^a B)=
\left\{
\begin{array}{ll}
q^{\frac{m-1}{2}}\frac{(-1+G_q(\eta'))}{2}, & \mbox{ if $\Tr_{q^m/q}(\omega^{a}b)\in C_0^{(2,q)}$, 
$\omega^a \in (H_1\cap \cdots \cap H_\frac{m-1}{2})\setminus A$, }\\
q^{\frac{m-1}{2}}\frac{(-1-G_q(\eta'))}{2}, & \mbox{ if $\Tr_{q^m/q}(\omega^{a}b)\in C_1^{(2,q)}$, 
$\omega^a \in (H_1\cap \cdots \cap H_\frac{m-1}{2})\setminus A$, 
}\\
\frac{q^{\frac{m-1}{2}}(q-1)}{2}, & \mbox{ if $\omega^a \in A \cap H_1\cap \cdots \cap H_\frac{m-1}{2}$,}
\\
0, & \mbox{ otherwise.}
 \end{array}
\right.
\end{align*}
\end{lemma}
\proof 
We compute the character values of $B$: 
\begin{align*}
\psi_{\F_{q^m}}(\omega^a B)=&\,\sum_{x_1,\ldots,x_{\frac{m-1}{2}}\in \F_q}
\sum_{y\in C_0^{(2,q)}}
\psi_{\F_{q^m}}(\omega^a a_1x_1)\cdots \psi_{\F_{q^m}}(\omega^a a_\frac{m-1}{2}x_\frac{m-1}{2})\psi_{\F_{q^m}}(\omega^a b y)\\
=&\,\left(\prod_{i=1}^{\frac{m-1}{2}}\sum_{x_i\in \F_q}
\psi_{\F_{q}}(\Tr_{q^m/q}(\omega^a a_i)x_i)\right)\left(\sum_{y\in C_0^{(2,q)}}\psi_{\F_q}(\Tr_{q^m/q}(\omega^a b) y)\right). 
\end{align*}
If  $\Tr_{q^m/q}(\omega^a a_i)\not=0$ for some $i=1,\ldots,\frac{m-1}{2}$, then it is clear that $\psi_{\F_{q^m}}(\omega^a B)=0$. Otherwise, we have 
\begin{align*}
\psi_{\F_{q^m}}(\omega^a B)=&\, q^{\frac{m-1}{2}}\sum_{y\in C_0^{(2,q)}}\psi_{\F_q}(\Tr_{q^m/q}(\omega^a b) y)\\
=&\, \left\{
\begin{array}{ll}
\frac{q^{\frac{m-1}{2}}(-1+G_q(\eta'))}{2}, & \mbox{ if $\Tr_{q^m/q}(\omega^{a}b)\in C_0^{(2,q)}$,  }\\
\frac{q^{\frac{m-1}{2}}(-1-G_q(\eta'))}{2}, & \mbox{ if $\Tr_{q^m/q}(\omega^{a}b)\in C_1^{(2,q)}$, }\\
\frac{q^{\frac{m-1}{2}}(q-1)}{2}, & \mbox{ if $\Tr_{q^m/q}(\omega^{a}b)=0$.}
 \end{array}
\right.
\end{align*}
Since $\Tr_{q^m/q}(\omega^a a_1)=\cdots =\Tr_{q^m/q}(\omega^a a_\frac{m-1}{2})=\Tr_{q^m/q}(\omega^a b)=0$ if and only if $\omega^a \in A \cap H_1\cap \cdots \cap H_\frac{m-1}{2}$, the assertion of the lemma follows. \qed

\vspace{0.3cm}
We are now ready to prove Theorem~\ref{thm:coSing}. 

{\bf Proof of Theorem~\ref{thm:coSing}:} \, 
From Lemmas~\ref{lem:coSing} and \ref{lem:coSing2}, we have 
\begin{align*}
\psi_{\F_{q^m}}(\omega^a\bigcup_{t\in \widehat{X}}C_t^{(2N,q^m)})=&\,\psi_{\F_{q^m}}(\omega^a D)=\sum_{i=1}^{\frac{m-1}{2}}\psi_{\F_{q^m}}(\omega^a T_\ell)
+\psi_{\F_{q^m}}(\omega^a B)\\
=&\, 
\left\{
\begin{array}{ll}
{\frac{(-1)^{i+\epsilon \frac{m-1}{2}}q^{\frac{m-1}{2}}(-1+(-1)^{j+\tau} G_q(\eta'))}{2},}, & \mbox{ if $\omega^a \in A_i\cap H_1\cap \cdots \cap H_{\ell-1}\cap H_{\ell,j}$ }\\
 & \mbox{ \hspace{1cm}for $\ell=1,\ldots,\frac{m-1}{2}$, $i,j=0,1$, }\\
\frac{q^{\frac{m-1}{2}}(-1+(-1)^i G_q(\eta'))}{2}, & \mbox{ if $\omega^a\in (H_1\cap \cdots \cap H_{\frac{m-1}{2}})\setminus A$ and }\\
 & \mbox{  \hspace{1cm}$\Tr_{q^m/q}(\omega^a b)\in C_i^{(2,q)}$ for $i=0,1,$}\\
0, & \mbox{ if $\omega^a \in A$. }
 \end{array}
\right.
\end{align*}
On the other hand, since $\psi_{\F_{q^m}}(\omega^a \bigcup_{t\in 2^{-1}I}C_t^{(N,q^m)})=\psi_{\F_{q^m}}(\omega^a D)+\psi_{\F_{q^m}}(\omega^{a+N} D)$, we have  
\[
\psi_{\F_{q^m}}(\omega^a \bigcup_{t\in 2^{-1}I}C_t^{(N,q^m)})\\
=
\left\{
\begin{array}{ll}
-(-1)^{i+\epsilon \frac{m-1}{2}}q^{\frac{m-1}{2}}, & \mbox{ if $\omega^a \in A_i\setminus (H_1\cap \cdots \cap H_{\frac{m-1}{2}})$, $i=0,1,$}\\
-q^{\frac{m-1}{2}}, & \mbox{ if $\omega^a\in (H_1\cap \cdots \cap H_{\frac{m-1}{2}})\setminus A$.}\\
0, & \mbox{ if $\omega^a \in A$.}
 \end{array}
\right.
\]
Hence, we have 
\begin{equation*}
2\psi_{\F_{q^m}}(\omega^a\bigcup_{t\in \widehat{X}}C_t^{(2N,q^m)})-\psi_{\F_{q^m}}(\omega^a \bigcup_{t\in 2^{-1}I}C_t^{(N,q^m)})=
\left\{
\begin{array}{ll}
0, & \mbox{ if $a\in 2^{-1}I\,(\mod{N})$,}\\
\pm q^{\frac{m-1}{2}}G_q(\eta'), & \mbox{ otherwise.}
 \end{array}
\right.
\end{equation*}
Thus, 
we conclude that  $\widehat{X}$ 
satisfies the condition \eqref{eq:3-chara2}  in Remark~\ref{rem:chara} (ii). 
\qed

\vspace{0.3cm}
As corollaries, we obtain the following. 
\begin{corollary}\label{cor:aff1}
For a prime power $q\equiv 3\,(\mod{4})$ and an odd integer $m>1$, there exists 
a $(q^{2m},r(q^{m}+1),-q^m+r^2+3 r,r^2+ r)$ strongly regular Cayley graph with $r=q^{m-1}(q-1)/2$. 
\end{corollary}
\proof
By Lemma~\ref{rem:quad-}, Remark~\ref{rem:chara}~(i)  and  Theorem~\ref{thm:coSing}, the corollary now follows. 
\qed
\begin{remark} {\em 
The strongly regular graph obtained in Corollary~\ref{cor:aff1} has the same parameter as the affine polar graph of elliptic type. Let $\Gamma $ 
be the strongly regular graph of Corollary~\ref{cor:aff1}  with $q=3$ and $m=3$. We checked by using a computer that  $\Gamma$ is {\bf not} isomorphic to the affine polar graph $\AP^-$ with the same parameters.  In particular, the size of the full automorphism group of $\Gamma$ (resp. $\AP^-$) is $2^2\cdot 3^7 \cdot 7$ (resp. $2^{10}\cdot 3^{12} \cdot 5\cdot 7$). }
\end{remark}

\begin{corollary}\label{cor:aff2}
For a prime power $q\equiv 1\,(\mod{4})$ and an odd integer $m>1$ such that $\gcd{(q-1,\frac{q^m-1}{q-1})}=1$, there exists a  $(q^{2m},r(q^m-1),q^m+r^2-3r,r^2-r)$ strongly regular Cayley graph, where  $r=q^{m-1}(q-1)/2$.
\end{corollary}
\proof
By Lemma~\ref{mainconstruction2+}, Remark~\ref{rem:chara}~(ii)  and  Theorem~\ref{thm:coSing}, the corollary now follows. 
\qed
\begin{remark} {\em 
The strongly regular graph obtained in Corollary~\ref{cor:aff2} has the same parameters as the affine polar graph of hyperbolic type. Let $\Gamma$ 
be the strongly regular graph of Corollary~\ref{cor:aff2}  with $q=5$ and $m=3$. We checked by using a computer that  $\Gamma$ is {\bf not} isomorphic to the affine polar graph $\AP^+$ with the same parameters. 
In particular, the size of the full automorphism group of $\Gamma$ (resp. $\AP^+$) is $2^3\cdot 5^6 \cdot 31$ (resp. $2^{11} \cdot 3^2\cdot 5^{12}\cdot 13 \cdot 31$). }
\end{remark}



\section*{Appendix: Proof of Lemma~\ref{lem:coSing}}
In this appendix, we give a proof of Lemma~\ref{lem:coSing}.  
\vspace{0.1in}

\noindent{\bf Proof of Lemma~\ref{lem:coSing}.} \, 
For $j=0$ or $1$, the characteristic function $g_{A,j}$ of $\{x\in \F_{q^m}\,|\,\Tr_{q^m/q}(x^2)\in C_j^{(2,q)}\}$ is given by 
\begin{equation}\label{eq:ch1}
g_{A,j}(x)=\frac{1}{q}\sum_{d\in \F_q}\sum_{s\in C_j^{(2,q)}}\psi_{\F_{q^m}}(dx^2)\psi_{\F_{q}}(-ds). 
\end{equation}
Similarly, the characteristic functions $g_{a_\ell}$ and $g_{a_\ell,j}$ of $\{x\,|\,\Tr_{q^m/q}(xa_\ell)=0\}$ and $\{x\,|\,\Tr_{q^m/q}(xa_\ell)\in C_j^{(2,q)}\}$ are,  respectively, given by 
\begin{equation}\label{eq:ch2}
g_{a_{\ell}}(x)=\frac{1}{q}\sum_{d\in \F_q}\psi_{\F_{q^m}}(dxa_\ell)
\end{equation}
and 
\begin{equation}\label{eq:ch3}
g_{a_{\ell},j}(x)=\frac{1}{q}\sum_{d\in \F_q}\sum_{s\in C_j^{(2,q)}}\psi_{\F_{q^m}}(dxa_\ell)\psi_{\F_{q}}(-ds), \, \quad  j=0,1. 
\end{equation}
We compute the character values $\psi_{\F_{q^m}}(\omega^a T_\ell)$. By the definition of $T_\ell$, we have 
\begin{equation}\label{eq:cha4}
\psi_{\F_{q^m}}(\omega^a T_\ell)=\sum_{j=0,1}\sum_{x\in \F_{q^m}}g_{A,j}(x)g_{a_1}(x)\cdots g_{a_{\ell-1}}(x)g_{a_{\ell},j}(x)\psi_{\F_{q^m}}(\omega^a x)
\end{equation}
By substituting \eqref{eq:ch1}, \eqref{eq:ch2} and \eqref{eq:ch3} into 
\eqref{eq:cha4}, we have 
\begin{equation}
\psi_{\F_{q^m}}(\omega^a T_\ell)
=\frac{1}{q^{\ell+1}}\sum_{x\in \F_{q^m}}\sum_{j=0,1}\sum_{d_0,d_1,\ldots,d_{\ell}\in \F_q}
\psi_{\F_{q^m}}(d_0x^{2}+(\omega^a+\sum_{i=1}^\ell d_ia_i)x)
\psi_{\F_q}(d_0 C_j^{(2,q)})\psi_{\F_q}(d_{\ell} C_j^{(2,q)}). \label{eq:compe}
\end{equation}
We compute the right hand side of \eqref{eq:compe} by dividing into the two  partial sums: $\Sigma_1$ and $\Sigma_2$, where $\Sigma_1$ is the contribution of the summands with $d_0=0$ and $\Sigma_{2}$ is the contribution of the summands  with $d_0\not=0$. 
Thus, $\psi_{\F_{q^m}}(\omega^a T_\ell)=\Sigma_1+\Sigma_{2}$. 

It is clear that 
\begin{align*} 
\Sigma_1=&\,\frac{q-1}{2q^{\ell+1}}\sum_{x\in \F_{q^m}}\sum_{j=0,1}\sum_{d_1,\ldots,d_{\ell}\in \F_q}
\psi_{\F_{q^m}}((\omega^a+\sum_{i=1}^\ell d_ia_i)x)
\psi_{\F_q}(d_{\ell} C_j^{(2,q)})\\
=&\,\left\{
\begin{array}{ll}
\frac{-q^{m-\ell -1}(q-1)}{2}, & \mbox{ if $\omega^a\in \langle a_1,\ldots,a_\ell\rangle\setminus \langle a_1,\ldots,a_{\ell-1}\rangle$,}\\
\frac{q^{m-\ell-1}(q-1)^2}{2}, & \mbox{ if $\omega^a\in \langle a_1,\ldots,a_{\ell-1}\rangle$,}\\
0, & \mbox{ otherwise.}
 \end{array}
\right.
\end{align*}
We next consider the partial sum  $\Sigma_{2}$. 
By Theorem~\ref{prop:charaadd},  
\begin{equation}\label{eq:outin}
\Sigma_{2}=\frac{G_{q^m}(\eta)}{q^{\ell+1}}\sum_{j=0,1}\sum_{d_0\in \F_q^\ast}\sum_{d_1,\ldots,d_{\ell}\in \F_q}
\psi_{\F_{q^m}}(-4^{-1}d_0^{-1}(\omega^a+\sum_{i=1}^\ell d_ia_i)^2)\eta(d_0)
\psi_{\F_q}(d_0 C_j^{(2,q)})
\psi_{\F_q}(d_{\ell} C_j^{(2,q)}), 
\end{equation}
where $\eta$ is the quadratic character of $\F_{q^m}$. 
Since $\Tr_{q^m/q}(a_{i}a_j)=0$ for $i,j\in \{1,\ldots,\ell\}$, we have 
\[
\Tr_{q^m/q}((\omega^a+\sum_{i=1}^\ell d_ia_i)^2)=\Tr_{q^m/q}(\omega^{2a}+ {2}\omega^a\sum_{i=1}^\ell d_ia_i). 
\]
Continuing from \eqref{eq:outin}, we have 
\[
\Sigma_{2}=\frac{G_{q^m}(\eta)}{q^{\ell+1}}\sum_{j=0,1}\sum_{h=0,1}(-1)^h\psi_{\F_q}(C_{j+h}^{(2,q)})\sum_{d_1,\ldots,d_{\ell}\in \F_q}
\psi_{\F_{q}}(-\Tr_{q^m/q}(\omega^{2a}+ {2}\omega^a\sum_{i=1}^\ell d_ia_i)C_h^{(2,q)})
\psi_{\F_q}(d_{\ell} C_j^{(2,q)}). 
\]
If  $\Tr_{q^m/q}(\omega^a a_i)\not=0$ for some $i=1,2,\ldots,\ell-1$, 
it is clear that $\Sigma_{2}=0$. Furthermore, if $\Tr_{q^m/q}(\omega^a a_\ell)=0$, it also holds that $\Sigma_{2}=0$. Thus, we assume that 
$\Tr_{q^m/q}(\omega^a a_i)=0$ for all $i =1,2,\ldots,\ell-1$ and 
$\Tr_{q^m/q}(\omega^a a_\ell)\not=0$. In this case, we have 
\begin{equation}\label{eq:P2}
\Sigma_{2}=\frac{G_{q^m}(\eta)}{q^2}\sum_{j=0,1}\sum_{h=0,1}(-1)^h\psi_{\F_q}(C_{j+h}^{(2,q)})\sum_{d_{\ell}\in \F_q}
\psi_{\F_{q}}(-\Tr_{q^m/q}(\omega^{2a}+ {2}\omega^ad_\ell a_\ell)C_h^{(2,q)})
\psi_{\F_q}(d_{\ell} C_j^{(2,q)}). 
\end{equation}
We compute the right hand side of \eqref{eq:P2} by dividing into the two partial sums: $\Sigma_{2.0}$ and $\Sigma_{2,1}$, where $\Sigma_{2,0}$ is the contribution of the summands with $d_\ell=0$ and $\Sigma_{2,1}$ is the contribution of the summands  with $d_\ell\not=0$. 
Thus, $\Sigma_2=\Sigma_{2,0}+\Sigma_{2,1}$. 

By \eqref{eq:Gaussquad}, we have 
\begin{align*}
\Sigma_{2,0}=&\,\frac{(q-1)G_{q^m}(\eta)}{2q^2}\sum_{j=0,1}\sum_{h=0,1}(-1)^h\psi_{\F_q}(C_{j+h}^{(2,q)})
\psi_{\F_{q}}(-\Tr_{q^m/q}(\omega^{2a})C_h^{(2,q)})\\
=&\, -\frac{(q-1)G_{q^m}(\eta)}{2q^2}\left\{
\begin{array}{ll}
G_q(\eta'), & \mbox{ if $-\Tr_{q^m/q}(\omega^{2a})\in C_0^{(2,q)}$,}\\
-G_q(\eta'), & \mbox{ if $-\Tr_{q^m/q}(\omega^{2a})\in C_1^{(2,q)}$,}\\
0, & \mbox{ if $\Tr_{q^m/q}(\omega^{2a})=0$,}
 \end{array}
\right.
\end{align*}
where $\eta'$ is the quadratic character of $\F_q$. On the other hand, 
by \eqref{eq:Gaussquad}, 
\begin{align*}
\Sigma_{2,1}
=&\,\frac{G_{q^m}(\eta)}{q^2}\sum_{j,h,k=0,1}(-1)^h
\psi_{\F_q}(C_{j+h}^{(2,q)})
\psi_{\F_{q}}(-\Tr_{q^m/q}(\omega^{2a}) {C_h^{(2,q)}})\psi_{\F_q}(-2\Tr_{q^m/q}(\omega^a a_\ell)C_{h+k}^{(2,q)})
\psi_{\F_q}(C_{j+k}^{(2,q)})\\
=&\,\frac{G_{q^m}(\eta)}{q^2}\left\{
\begin{array}{ll}
\frac{G_{q}(\eta')(-1+G_q(\eta')^3)}{2}, & \mbox{ if $-\Tr_{q^m/q}(\omega^{2a})\in C_0^{(2,q)}$, $-2\Tr_{q^m/q}(\omega^aa_\ell)\in C_0^{(2,q)}$,}\\
\frac{-G_{q}(\eta')(1+G_q(\eta')^3)}{2}, & \mbox{ if $-\Tr_{q^m/q}(\omega^{2a})\in C_0^{(2,q)}$, $-2\Tr_{q^m/q}(\omega^a a_\ell)\in C_1^{(2,q)}$,}\\
\frac{-G_{q}(\eta')(-1+G_q(\eta')^3)}{2}, & \mbox{ if $- \Tr_{q^m/q}(\omega^{2a})\in C_1^{(2,q)}$, $-2\Tr_{q^m/q}(\omega^aa_\ell)\in C_0^{(2,q)}$,}\\
\frac{G_{q}(\eta')(1+G_q(\eta')^3)}{2}, & \mbox{ if $-\Tr_{q^m/q}(\omega^{2a})\in C_1^{(2,q)}$, $-2\Tr_{q^m/q}(\omega^a a_\ell)\in C_1^{(2,q)}$,}\\
0, & \mbox{ if $\Tr_{q^m/q}(\omega^{2a})=0$.}
 \end{array}
\right.
\end{align*}
Noting that $G_{q^m}(\eta)=G_q(\eta')^m$ and $G_{q}(\eta')^2=(-1)^\epsilon q$, we have 
\begin{align*}
\Sigma_2=&\,\Sigma_{2,0}+\Sigma_{2,1}\\
=&\,(-1)^{\epsilon \frac{m+1}{2}}q^{\frac{m-1}{2}}\left\{
\begin{array}{ll}
\frac{-1+(-1)^\epsilon G_q(\eta')}{2}, & \mbox{ if $-\Tr_{q^m/q}(\omega^{2a})\in C_0^{(2,q)}$,  $-2\Tr_{q^m/q}(\omega^a a_\ell)\in C_0^{(2,q)}$,}\\
\frac{-1-(-1)^\epsilon G_q(\eta')}{2}, & \mbox{ if $-\Tr_{q^m/q}(\omega^{2a})\in C_0^{(2,q)}$, $-2\Tr_{q^m/q}(\omega^aa_\ell)\in C_1^{(2,q)}$,}\\
\frac{1-(-1)^\epsilon G_q(\eta')}{2}, & \mbox{ if $-\Tr_{q^m/q}(\omega^{2a})\in C_1^{(2,q)}$,  $-2\Tr_{q^m/q}(\omega^aa_\ell)\in C_0^{(2,q)}$,}\\
\frac{1+(-1)^\epsilon G_q(\eta')}{2}, & \mbox{ if $-\Tr_{q^m/q}(\omega^{2a})\in C_1^{(2,q)}$,  $-2\Tr_{q^m/q}(\omega^{a}a_\ell)\in C_1^{(2,q)}$,}\\
0, & \mbox{ if $\Tr_{q^m/q}(\omega^{2a})=0$.}
 \end{array}
\right.
\end{align*}
This completes the proof of the lemma. \qed
\end{document}